\documentclass[11pt, amsfonts,amssymb, amsmath]{amsart}
\usepackage{times}
\usepackage{amsfonts}
\usepackage{amssymb} 
\theoremstyle{plain}
\newtheorem{theorem}{Theorem}[section]
\newtheorem{corollary}[theorem]{Corollary}
\newtheorem{lemma}[theorem]{Lemma}

\theoremstyle{definition}
\newtheorem{definition}[theorem]{Definition}
\newtheorem{remark}[theorem]{Remark}
\newtheorem{hypothesis}[theorem]{Hypothesis}

\newenvironment{nolabels}
{\begin{list}{}{\leftmargin 15pt}}
{\end{list}}

\newcommand{\sq}[2]{\sideset{^{#1}}{}{\operatorname{#2}}}
\newcommand{\md}{\operatorname{Mod}}
\newcommand{\ec}{\operatorname{EC}}
\newcommand{\pc}{\operatorname{PC}}
\newcommand{\tp}{\operatorname{tp}}
\newcommand{\rang}{\operatorname{rang}}
\newcommand{\dom}{\operatorname{dom}}
\newcommand{\aut}[1]{\operatorname{Aut}_{#1}(\mathfrak{C})}
\newcommand{\av}{\operatorname{Av_D}}
\newcommand{\fml}{\operatorname{Fml}}
\newcommand{\conc}{\Hat{\ }}

\newcommand{\union}{\cup}
\newcommand{\intersect}{\cap}
\newcommand{\Union}{\bigcup}

\newcommand{\al}[1]{{\aleph}_{#1}}

\newcommand{\infinity}{\infty}

\newcommand{\bpf}{\begin{proof}}
\newcommand{\epf}{\end{proof}}
\renewcommand{\empty}{\emptyset}

\newbox\noforkbox \newdimen\forklinewidth
\setbox0\hbox{$\textstyle\smile$}
\setbox1\hbox to \wd0{\hfil\vrule width \forklinewidth depth-2pt
 height 10pt \hfil}
\setbox\noforkbox\hbox{\lower 2pt\box1\lower 2pt\box0\relax}
\def\unionstick{\mathop{\copy\noforkbox}\limits}

\def\nonfork_#1{\unionstick_{\textstyle #1}}

\setbox0\hbox{$\textstyle\smile$}
\setbox1\hbox to \wd0{\hfil{\sl /\/}\hfil}
\setbox2\hbox to \wd0{\hfil\vrule height 10pt depth -2pt width
               \forklinewidth\hfil}
\newbox\doesforkbox
\setbox\doesforkbox\hbox{\lower 2pt\box1 \lower 2pt\box2\lower2pt\box0\relax}
\def\nunionstick{\mathop{\copy\doesforkbox}\limits}

\def\fork_#1{\nunionstick_{\textstyle #1}}

\begin{document}

\title{Ranks and Pregeometries in
Finite Diagrams}

\author{Olivier Lessmann} 
\address{Department of Mathematical Sciences \\
Carnegie Mellon University \\
Pittsburgh, PA 15213}
\email{lessmann@andrew.cmu.edu}
\date{\today}

\thanks{This work is a chapter of the author's PhD thesis,
under Prof. Rami Grossberg. I am deeply grateful to him 
for his constant guidance and support.}

\begin{abstract} The study of classes of models 
of a finite
diagram was initiated by S. Shelah in
1969. A \emph{diagram} $D$ is a set of types over the empty set,
and the class of models of the diagram $D$ 
consists of the models of $T$ which omit all the
types not in $D$. 
In this work, we introduce 
a natural dependence relation on the subsets of the models for 
the $\aleph_0$-stable case
which share many of the formal
properties of forking. 
This is achieved by considering a rank for this framework
which is bounded when the diagram $D$ is  $\aleph_0$-stable.
We can also obtain pregeometries
with respect to this dependence relation.
The dependence relation is the natural one induced by
the rank, and the pregeometries exist on the set
of realizations of types of minimal rank.
Finally, these concepts are used to generalize many of the classical
results for models of a totally transcendental first-order theory. 
In fact, strong analogies arise:
 models are determined by their pregeometries or their
relationship with their pregeometries; however
 the proofs are different, as we do not have
compactness.
This is illustrated with positive results 
(categoricity) as well as negative results
(construction  of nonisomorphic models). 
\end{abstract}

\maketitle

\section{Introduction}

The problem of categoricity has been a driving force
in model theory since its early development in the late 1950's.
For the countable first-order case, M.~Morley
in 1965 (\cite{Mo}) introduced a rank which captures $\aleph_0$-stability,
and used it to construct prime models and give a proof of
{\L}o\v{s} conjecture. 
In 1971, J.~Baldwin and A.~Lachlan \cite{BlLa} gave an
alternative proof using the fact that algebraic closure
induces a pregeometry on strongly minimal sets.
Their proof generalizes ideas from Steinitz's
famous 1910 theorem of categoricity for
algebraically closed fields.
{\L}o\v{s} conjecture for uncountable
languages was solved in 1970 by S.~Shelah \cite{Sh:70} introducing
a rank which corresponds to the superstable case. 
Later, Shelah discovered a dependence relation
called forking and more general pregeometries, 
and since then, these ideas have been extended
to more and more general first-order contexts, each of them
corresponding to a specific rank: $\aleph_0$-stable, 
superstable, stable and simple.

The problem of categoricity for non-elementary classes 
is quite considerably more involved. 
In 1971, H.~J.~Keisler (see \cite{Ke}) 
proved a categoricity theorem for Scott sentences 
$\psi \in L_{\omega_1\omega}$, which in a sense 
generalizes Morley's Theorem. 
To achieve this, Keisler made the additional assumption that 
$\psi$ admits $\aleph_1$-homogeneous models. 
Later, L.~Marcus, with the assistance of Shelah (see \cite{MaSh}),
produced an example of a categorical 
$\psi \in L_{\omega_1\omega}$ that does not have any
$\aleph_1$-homogeneous model, so this is not the most general case.
Since then, many of Shelah's hardest papers in model theory
have been dedicated the
categoricity problem and to 
the development of general classification theory for
non-elementary classes.
Among the landmarks, one should mention \cite{Sh:4} about sentences in 
$L_{\omega_1\omega}(Q)$ 
which answers a question of Harvey Friedman's list 
(see \cite{Fr}).
In \cite{Sh:871} and \cite{Sh:872} 
a version of Morley's Theorem is proved for a special kind of 
formulas $\psi \in L_{\omega_1\omega}$ which are called excellent. 
It is noteworthy that to deal with these non-elementary classes, 
these papers introduced several crucial ideas, among
them stable amalgamation, $2$-goodness and others, 
which are now essential parts
of the proof of the ``Main Gap" for first-order, countable theories. 
Later, R.~Grossberg and B.~Hart completed the classification of 
excellent classes
and 
gave a proof of the Main Gap for those classes (\cite{GrHa}).
H.~Kierstead also continued the study 
of sentences in $L_{\omega_1\omega} (Q)$ (see \cite{Ki}). 
He introduced
a generalization of strongly minimal formulas by replacing 
``non-algebraic'' by ``there exists uncountably many" 
and obtained results about countable models of 
these classes using \cite{Sh:4}.
In \cite{Sh:300}, Shelah began the classification theory for universal classes 
(see also ICM 1986/videotape) and 
is currently  working on a book entirely dedicated to them.
He also started the classification of classes in a context somewhat 
more general than $\pc(T_1, T, \Gamma)$,
see \cite{Sh:88}, \cite{Sh:576} and \cite{Sh:600}.
In a related work,
Grossberg started studying the
classification of $\md(\psi)$ for $\psi \in L_{\lambda^+\omega}$ under
the assumption that there exists a ``Universal Model" for $\psi$ and
studied relatively saturated substructures (see \cite{Gr:1} and \cite{Gr:2}). 
This seems to be a natural hypothesis which others have made as well
(for example \cite{Sh:88},
\cite{KlSh} and \cite{BlSh:3}).
As a matter of fact, it is conjectured that if an abstract
class of models $\mathcal{K}$ is categorical
above the Hanf number,
then $\mathcal{K}$ has the $\mu$-amalgamation property for every $\mu$  
(this implies the existence of $\mu^+$-universal models, under
the General Continuum Hypothesis).

There are several striking differences between the problem
of categoricity for first-order and the non-elementary case.
First, it appears that
classification for non-elementary classes is
sensitive to the axioms of set theory.
Second, the methods used are heavily combinatorial:
there is no ``forking'' (though splitting and
strong splitting are sometimes well-behaved), 
and the use of pregeometries to understand systematically
models of a given class is virtually absent. (A nice example
of pregeometries is hidden in the last
section of \cite{Sh:4} and only \cite{Ki} 
has used them to study countable models.)
However, stability was not developed 
originally for first-order.
In 1970, Shelah published \cite{Sh:1}, 
where he introduced some of the most 
fundamental ideas of classification theory (stability, splitting of types, 
existence of indiscernibles, several notions of prime models etc.). 
In this paper, Shelah considered classes of models which omit all types in 
$D(T)-D$, for a fixed \emph{diagram} $D \subseteq D(T)$. 
This class is usually denoted $\ec(T, \Gamma)$, where $\Gamma$ stands
for $D(T)-D$. 
He made assumptions of two kinds 
(explicitly in his definition of stability): 
(1) restriction on the cardinality
of the space of types realizable by the models, and (2) 
existence of models realizing many
types. In fact, the context studied by Keisler
in his categoricity
result for $L_{\omega_1 \omega}$, turns out to be
 the $\aleph_0$-stable case in the above sense.
This is made precise by the following results. 
(C.-C.~Chang:) The class of models
of a sentence $\psi \in L_{\omega_1\omega}$ is equal to 
the class $\pc(T_1, T, \Gamma)$, which is the class of reducts to $L(T)$ 
of models of a first-order countable theory $T_1$ containing $T$, and
 omitting a set of types 
$\Gamma \subseteq D(T_1)$. 
(Shelah:) The number of models of a Scott sentence 
$\psi \in L_{\omega_1\omega}$ is equal to the
number of models of $\ec(T, \Gamma)$, for some countable $T$,
where $\Gamma$ the set of isolated types of $T$.

In  retrospect, it seems that what prevented the emergence
of a smooth theory for $\aleph_0$-stable diagrams
is the absence of a rank like Morley's rank.
Considering the success of the use of pregeometries
to understand models in the first-order 
$\aleph_0$-stable case, if one hopes to lift these
ideas to more general contexts, it appears
that $\aleph_0$-stable diagrams constitute
a natural test case.
This is the main goal of this paper.
 We try to develop 
what Shelah calls the structure part of the theory for
the class $\ec(T, \Gamma)$, under the assumption
that it is $\aleph_0$-stable (in the sense of \cite{Sh:1}).
In fact, as in \cite{Sh:2}, we assume that $\ec(T,\Gamma)$
contains a large homogeneous model 
(which follows from Shelah's original definition
of stability for $\ec(T, \Gamma)$, see Theorem 3.4. in \cite{Sh:1}), so that
the stability assumptions only deal with the cardinality of 
the spaces of types.
This hypothesis allows us to do all the work in ZFC, 
in contrast to \cite{Sh:4}, 
\cite{Sh:871}, \cite{Sh:872} or \cite{Ki} for example.

The paper is organized as follows. 
\begin{description}
\item[Section 1]
We describe the general context.
\item[Section 2] 
We introduce a rank for this framework
which captures $\aleph_0$-stability 
(it does not generalize Morley rank, 
but rather generalizes what Shelah calls $R[p, L, 2]$). 
This rank differs from previously studied ranks in two ways:
(1) it allows
us to deal with general diagrams (as opposed to the atomic case or 
the first-order case)
and (2) the definition is relativized to a given set (which allows
us to construct prime models).
By analogy with the first-order case, 
we call $D$ totally transcendental when the rank is
bounded. 
For the rest of the paper, we only consider totally 
transcendental $D$, and we make no assumption
on the cardinality of $T$.
We study the basic properties of this rank, and examine the natural
 dependence relation that it induces on the subsets of the models.
We are then able to obtain many of the classical properties of forking, which
we summarize in Theorem \ref{forkingproperties}.
We also obtain stationary types with respect 
to this dependence relation, and they
turn out to behave well: they satisfy in addition the symmetry property,  
and can be represented by averages. 
\item[Section 3]
We focus on pregeometries. 
Regular types are defined in the usual manner 
(but with this dependence relation 
instead of forking, of course), and the dependence relation 
on the set of realizations
of a regular type yields a pregeometry. 
We can show that stationary types of minimal rank
are regular, and this is used to show that they exist very often.
We also consider a more concrete kind of regular types, 
which are called minimal. 
They could be defined independently 
by replacing ``non-algebraic" by ``realized outside any model
which contains the set of parameters" in the usual definition of 
strongly minimal
formulas. (This can be done for any suitable class of models, 
as in the last section
of~\cite{Sh:4}.) We could show directly that the natural closure 
operator induces a pregeometry on the set of realizations 
in any $(D,\aleph_0)$-homogeneous model.
We choose not to do this, 
and instead we consider minimal types only when the natural
dependence relation coincides with the one given by the rank.
This allows us to use the results we have already obtained and have a picture
which is conceptually similar to the first-order totally transcendental case
(where strongly minimal types are stationary and regular, and 
the unique nonforking extension is also the unique non-algebraic one).
Another reason is that the proofs are identical to those
which use the rank, and this
presentation permits us to skip them. 
\item[Section 4]
Here, we give various applications 
of both the rank and the pregeometries to the class 
$\mathcal{K}$ of $(D,\aleph_0)$-homogeneous 
models of a totally transcendental diagram. 
We introduce unidimensionality for diagrams.
We are able to adapt techniques of
Baldwin-Lachlan (see \cite{BlLa}) to our context for the categoricity proof.
In fact, we obtain a picture strikingly similar to 
the first-order totally transcendental
case. (1) If $D$ is totally transcendental,
then over any $D$-set there is a prime model for $\mathcal{K}$
(this improves parts of Theorems 5.3 and 5.10 of \cite{Sh:1}).
(2)~If $D$ is totally transcendental, then
$\mathcal{K}$ is categorical in some $\lambda > |T|+|D|$
if and only if $\mathcal{K}$ is categorical in every $\lambda > |T|+|D|$ 
if and only if every model of $\mathcal{K}$ is
prime and minimal over the set of realizations of a minimal type 
if and only if
every model of $\mathcal{K}$  of cardinality $>|T|+|D|$ is $D$-homogeneous.
(3) If $D$ is totally transcendental and if 
there is a model of $\mathcal{K}$ of cardinality above $|T|+|D|$
which is not $D$-homogeneous, then for
any $|T|+|D| \leq \mu \leq \lambda$, there exists maximally 
$(D,\mu)$-homogeneous
models in $\mathcal{K}$ of cardinality $\lambda$ (see the definition below).
If $T$ is countable this implies, in particular, that
for each ordinal $\alpha$ 
the class
$\mathcal{K}$ has at least $|\alpha|$ models of cardinality $\aleph_\alpha$.
When $|T| < 2^{\aleph_0}$, the categoricity assumption on $\mathcal{K}$ implies
that $D$ is totally transcendental, if $D$ is
the set of isolated types of $T$. As a byproduct, this gives an alternative
proof to Keisler's theorem
which  works so long as $|T|< 2^{\aleph_0}$ (whereas Keisler's soft
$L_{\omega_1 \omega}$ methods do not generalize).
\end{description}
Using regular types and prime models,
we could also give a decomposition theorem,
 but we do not include it here since it is a particular case of
a more general abstract decomposition theorem, part of a joint
work with R. Grossberg.

\section{The Context}

Let $T$ be a first-order theory in the language $L(T)$.
Let $\overline{M}$ be a very large saturated model of $T$. 
All sets are assumed to be subsets of $\overline{M}$.
As usual, 
\[
\tp(\bar{c},A)=\{\, \phi(\bar{x}, \bar{a}) \mid
 \overline{M} \models \phi[\bar{c},\bar{a}], 
\ell(\bar{c})=\ell(\bar{x}), \phi \in L(T) \,\}.
\]  
We say that $p(\bar{x})$ is a \emph{complete type over $A$} in $n$ variables 
if $\ell(\bar{x})=n$ and there is 
$\bar{c}$ in $\overline{M}$ such that 
$p(\bar{x})=\tp(\bar{c},A)$. 
The \emph{ diagram of $T$}, denoted by $D(T)$, is the set of complete 
types over
the empty set.
$S^n(A)$ is the set of all complete types over $A$ in $n$ variables.
$S^1(A)$ is written $S(A)$. 
Given a set of formulas $p$, we let $\dom(p)$ 
be the set of parameters appearing
in the formulas of $p$. 
We say that $p$ is over $A$ if $\dom(p)$ is contained in $A$.
Finally, given a type $p$ and a model $M$, we denote by
$p(M)$ the set of realizations of $p$ in $M$.

The following notions of \emph{diagram $D$} were defined by Shelah
in \cite{Sh:1}.

\begin{definition}

\begin{enumerate}
\item 
For any set $A$, let 
$D(A)=\{\, \tp(\bar{c},\empty)
\mid \bar{c} \in A \,\} \subseteq D(T)$;

\item 
For a model $M$ of $T$, let $D(M)=D(|M|)$.
\end{enumerate}

\end{definition}

\begin{definition}\label{dmodel}  Let $D \subseteq D(T)$.

\begin{enumerate}

\item 
$A$ is called a \emph{$D$-set} if $D(A) \subseteq D$;

\item 
A model $M$ of $T$ is called a \emph{$D$-model} if $D(M) \subseteq D$;

\item 
Define 
$S_D(A)=\{\, p \in S(A) \mid \text{ if $\bar{c} \models p$ then 
$A \cup \bar{c}$ is a $D$-set }\,\}$.

\end{enumerate}

\end{definition}

\begin{remark} $|S_D(A)|=|S^n_D(A)|$ provided both are infinite, 
so we will usually not write
the superscript.
\end{remark}

Here, we follow \cite{Sh:2}.

\begin{definition} Let $D \subseteq D(T)$.\label{hom}

\begin{enumerate}

\item 
The diagram $D$ is called \emph{stable in $\lambda$} if for any 
$D$-set $A$ of cardinality at most
$\lambda$, we have $|S_D(A)|\leq \lambda$;

\item 
The diagram $D$ is called \emph{stable} if there is $\lambda$ such that 
$D$ is stable in $\lambda$,
and we say that $D$ is \emph{unstable} if $D$ is not stable;

\item \label{homog} 
A $D$-model $M$ is called \emph{$(D,\lambda)$-homogeneous} if $M$ realizes
every type $p \in S_D(A)$ over subsets $A$ of $|M|$ 
of cardinality less than $\lambda$;

\item 
A $D$-model $M$ is \emph{$D$-homogeneous} if $M$ is $(D,\|M\|)$-homogeneous.

\end{enumerate}
\end{definition}

The following definition is due to 
Grossberg and Shelah in \cite{GrSh:2}.

\begin{definition}  We say that $D$ has the \emph{$\infinity$-order property} 
if
 for every $\lambda$, there is a formula
$\phi(\bar{x}, \bar{y}, \bar{z})$, a sequence $\bar{c}$ and
a set of sequences  $I=\{\, \bar{a}_i \mid  i < \lambda\,\}$, 
such that the following two conditions hold:
\begin{enumerate}
\item
$I \cup \bar{c}$ is a $D$-set;
\item
$\models \phi[\bar{a}_i, \bar{a}_j, \bar{c}]$ if and only if 
$i < j < \lambda$.
\end{enumerate}
\end{definition}

\begin{theorem} \cite{GrSh:2}\label{order} 
$D$ has the $\infinity$-order property if and only if 
there is a formula $\phi(\bar{x}, \bar{y}, \bar{z})$,
a sequence $\bar{c}$ and a set of sequences 
$I=\{\, \bar{a}_i \mid i < \beth_{(2^{|T|})^+} \,\}$, 
such that the following two conditions hold:
\begin{enumerate}
\item
$I \cup \bar{c}$ is a $D$-set;
\item
$\models \phi[\bar{a}_i, \bar{a}_j, \bar{c}] 
\text{ if and only if }
i < j <\beth_{(2^{|T|})^+}$.
\end{enumerate}
\end{theorem}

\begin{definition} Let $D \subseteq D(T)$ and let $\Gamma= D(T)-D$.
Define
\[
\ec(T,\Gamma)=\{\, M \models T \mid M \text{ omits every type in } \Gamma \,\}.
\]
Equivalently, 
\[
\ec(T,\Gamma)=\{\, M \models T \mid M \text{ is a $D$-model }\,\}.
\]
\end{definition}

For the rest of the paper, we will study 
the class $\ec(T, \Gamma)$, where $\Gamma=D(T)-D$ for a fixed diagram
$D \subseteq D(T)$,
under the following hypothesis.

\begin{hypothesis} There exists a $(D, \chi)$-homogeneous
model $\mathfrak{C} \in \ec(T,\Gamma)$ for some $\chi$ larger 
than any cardinality mentioned in this paper.
\end{hypothesis}
 
This implies that all $D$-models can be assumed to sit inside $\mathfrak{C}$,
and that model satisfaction is with respect to $\mathfrak{C}$.
In this context, Shelah proved the following results.

\begin{theorem} [The Stability Spectrum]\cite{Sh:1} \label{stabs} 
One of the following conditions must hold:
\begin{enumerate}
\item 
$D$ is unstable;
\item
 There are 
$\kappa(D) \leq \lambda(D) < \beth_{(2^{|T|})^+}$ such that for every $\mu$,
$D$ is stable in $\mu$ if and only if $\mu \geq \lambda(D)$ and 
$\mu^{<\kappa(D)}=\mu.$
\end{enumerate}
\end{theorem}

\begin{theorem} [The Homogeneity Spectrum]\label{homogeneity} \cite{Sh:2}\\
There is a $D$-homogeneous model of cardinality $\lambda$
if and only if
$\lambda \geq |D|$ and $D$ is stable in  $\lambda$ or 
$\lambda^{<\lambda}=\lambda.$
\end{theorem}

For an alternative and self-contained exposition of above two theorems, 
see \cite{GrLe}.

In the same paper, Shelah proved the following theorem. 
We will make use of a particular case which we will prove using the
rank.

\begin{theorem} \cite{Sh:2} \label{union} 
Let $D$ be stable.
If $\langle M_i \mid i < \alpha \rangle$ is an increasing sequence of
$(D,\mu)$-homogeneous models and
the cofinality of $\alpha$ is at least $\kappa(D)$, 
then $\Union_{i<\alpha} M_i$ 
is $(D,\mu)$-homogeneous.
\end{theorem}

The next theorem will be used to show the symmetry property
of the rank.

\begin{theorem} \label{unstable} \cite{Sh:5}
$D$ is unstable if and only if $D$ has the $\infinity$-order property.
\end{theorem} 

Using \cite{Sh:1} together with the method of \cite{Sh:a} Theorem 2.12 
and Theorem \ref{unstable}, one can easily show:

\begin{theorem} 
If $D$ is stable in $\lambda$, 
$A$ is a $D$-set of cardinality at most $\lambda$,
and $I$ is a $D$-set of finite 
sequences of cardinality at least $\lambda^+$,
then there is $J \subseteq I$ of cardinality $\lambda^+$,
 such that $J$ is an indiscernible set over $A$.
\end{theorem}

We will use the following properties of $\kappa(D)$ in the case when
$\kappa(D)=\aleph_0$, and we will actually provide alternative
proofs to these facts using the rank.

\begin{definition} Suppose $D$ is stable, $I$ is a $D$-set, 
which is a set of indiscernibles and $A$ is a $D$-set.
Define
\[
\av(I,A)= \{\,
\phi(\bar{x}, \bar{a})\mid 
\bar{a} \in A, \phi(\bar{x},\bar{y}) \in L(T) 
\text{ and }
|\phi(I, \bar{a})| \geq \kappa(D) \,\}.
\]
\end{definition}

\begin{lemma} \cite{Sh:2}\label{kappa} 
Suppose $D$ is stable, $I$ is a $D$-set, 
which is a set of indiscernibles and $A$ is a $D$-set.
Then
\begin{enumerate}
\item 
$\av(I,A) \in S_D(A)$;

\item
There exists $J$ a subset of $I$ with $|J|< |A|^++\kappa(D)$ such
that $I-J$ is indiscernible over $A \cup J$;

\item
If $|I| \geq |A|^+ + \kappa(D)$, then there is $\bar{a}$ in $I$
realizing $\av(I,A)$.

\end{enumerate}
\end{lemma}

\section{Rank, Stationary Types and Dependence relation}

We first introduce a rank for the class of $D$-models
 (see Definition \ref{dmodel}) 
which generalizes the rank from \cite{Sh:871}.  
We then prove basic properties of it
which show that it is well-behaved and is natural for this class. 

\begin{definition} \label{rank}
For any set of formulas $p(\bar{x},\bar{b})$ with parameters in 
$\bar{b}$, and $A$ a subset of $\mathfrak{C}$ containing $\bar{b}$,
 we define the \emph{rank $R_A[p]$}.
The rank  $R_A[p]$ will be an ordinal, $-1$, or
$\infinity$ and we have the usual ordering
$-1 < \alpha < \infinity$ for any ordinal $\alpha$.
 We define the relation $R_A[p] \geq \alpha$ 
by induction on  $\alpha$.

\begin{enumerate}

\item 
$R_A[p] \geq 0$ 
if $p(\bar{x},\bar{b})$ is realized in 
$\mathfrak{C}$;

\item 
$R_A[p] \geq \delta$, when $\delta$ is a limit ordinal,
if $R_A[p] \geq \alpha$ for every $\alpha < \delta$;

\item  
$R_A[p] \geq \alpha + 1$  if the following two conditions hold:

\begin{enumerate}

\item 
There is $\bar{a} \in A$ and a formula $\phi(\bar{x},\bar{y})$
such that 
\[
R_A[p \cup \phi(\bar{x},\bar{a})] \geq \alpha 
\quad 
\text{ and }
\quad
R_A[p \cup \neg \phi(\bar{x},\bar{a})] \geq \alpha;
\]
\item 
For every $\bar{a} \in A$ there is $q(\bar{x},\bar{y}) \in D$ 
such that 
\[
R_A[p \cup q(\bar{x},\bar{a})] \geq \alpha.
\]

\end{enumerate}

\end{enumerate}
We write:
\begin{nolabels}
\item
$R_A[p]=-1$ if $p$ is not realized in $\mathfrak{C}$;
\item
$R_A[p]=\alpha$ if $R_A[p] \geq \alpha$ but it is not
the case that $R_A[p] \geq \alpha + 1$;
\item
$R_A[p]=\infinity$ 
if $R_A[p] \geq \alpha$ for every ordinal $\alpha$. 

\end{nolabels}

For any set of formulas $p(\bar{x})$ over 
$A \subseteq \mathfrak{C}$, we let 
\[
R_A[p]=\min \{\, R_A[q] \mid 
q \subseteq p\restriction B, B \subseteq \dom(p), \text{ $B$ finite }\,\}.
\]
We omit the subscript $A$ when $A = \mathfrak{C}$.
\end{definition}

We need several basic properties of this rank. 
Some of them are purely technical and are stated
here for future reference.
Most of them are analogs of the usual properties for ranks in the
first-order case, with the exception of (2)
and (3).  The proofs 
vary from the first-order context because of the
second clause at successor stage, but they are 
all routine inductions. 

\begin{lemma} \label{rankproperties} 
Let $A$ be a subset of  $\mathfrak{C}$.
\begin{enumerate}

\item 
$R_A[\{\,\bar{x}=\bar{c}\,\}]=0$.

\item \label{goodness} 
If $p$ is over a finite set or $p$ is complete, 
then $R_A[p]\geq 0$ if and only if there is $B \subseteq A$ and $q \in S_D(B)$
 such that $p \subseteq q$.

\item \label{2}
If $A$ is $(D,\aleph_0)$-homogeneous and  
$\tp(\bar{a}, \empty)=\tp(\bar{b}, \empty)$
(for $\bar{a}, \bar{b}\in A$), then
$R_A[p(\bar{x},\bar{b})]=R_A[p(\bar{x},\bar{a})]$.

\item \textup{(Monotonicity)}\label{implication} 
If $p \vdash q$ and $p$ is over a finite set, 
then $R_A[p] \leq R_A[q]$. 

\item \label{rankinvariance} 
If $p$ is over $B \subseteq A$ and $f \in \aut{}$
then $R_A[p]=R_{f(A)}[f(p)]$.

\item \textup{(Monotonicity)} \label{subset} 
If $p \subseteq q$ then $R_A[p] \geq R_A[q]$.

\item \textup{(Finite Character)}\label{finite}
There is a finite $B \subseteq  \dom(p)$
such that 
\[
R_A[p]=R_A[p \restriction B].
\]

\item \label{8} 
If $R_A[p]= \alpha$ and $\beta < \alpha$, 
then there is $q$ over $A$ such that $R_A[q]=\beta$.

\item \label{bound} 
If $R_A[p] \geq (|A| +2^{|T|})^+$, then $R_A[p]=\infinity$. \\
Moreover, when $A$ is $(D,\aleph_0)$-homogeneous, the bound is $(2^{|T|})^+$.

\end{enumerate}
\end{lemma}

\bpf (1) Trivial

(2) 
Suppose $p \subseteq q \in S_D(B)$, and $B \subseteq A$.
Since $\mathfrak{C}$ is $(D, \chi)$-homogeneous,
and $q \in S_D(B)$, then $q$ is realized in $\mathfrak{C}$.
Hence $p$ is realized in $\mathfrak{C}$ and $R_A[p] \geq 0$.

For the converse, if $p$ is over a finite set, and $R_A[p]\geq 0$, 
then there is $\bar{c} \in \mathfrak{C}$ 
realizing $p$. 
Thus $\tp(\bar{c}, \dom(p))$ extends $p$ and
$\tp(\bar{c}, \dom(p)) \in S_D(\dom(p))$.

If $p$ is complete, then
there is $B \subseteq A$ such that $p \in S(B)$. 
Now let $\bar{c}$ (not necessarily in $\mathfrak{C}$) realize
$p$. 
For every $\bar{b} \in B$,  $R_A[p\restriction \bar{b}] \geq 0$, and 
so there is $\bar{c}' \in \mathfrak{C}$ 
realizing $p \restriction \bar{b}$.
But $\tp(\bar{c},\bar{b})=
p \restriction \bar{b}=\tp(\bar{c}',\bar{b})$ 
since $p$ is complete. 
Thus $\tp(\bar{c}\bar{b}, \empty) \in D$, so $p \in S_D(B)$.

(3) 
By symmetry, it is enough to show that for every ordinal $\alpha$, 
\[
 R_A[p(\bar{x},\bar{b})]\geq \alpha 
\qquad
\text{ implies }
\qquad
R_A[p(\bar{x},\bar{a})]\geq \alpha.
\] 
We prove that this is true for all types
 by induction on $\alpha$.
\begin{itemize}
\item 
When $\alpha =0$, we know that there is $\bar{c} \in \mathfrak{C}$ 
realizing 
$p(\bar{x},\bar{a})$. 
Then, since $\tp(\bar{a}, \empty)=\tp(\bar{b}, \empty)$ and 
$A$ is  $(D,\aleph_0)$-homogeneous, there is $\bar{d} \in A$ such
that $\tp(\bar{c}\bar{a}, \empty)=
\tp(\bar{d}\bar{b}, \empty)$. 
But then 
$p(\bar{x},\bar{b}) \subseteq \tp(\bar{d},\bar{b})$. 
Hence $p(\bar{x},\bar{b})$ is realized in $\mathfrak{C}$, 
so $R_A[p(\bar{x},\bar{b})] \geq 0$.
\item 
When $\alpha$ is a limit ordinal, this is true by induction.
\item 
Suppose $R_A[p(\bar{x},\bar{a})]\geq \alpha + 1$.
First, there is 
$\bar{c} \in A$ and $\phi(\bar{x},\bar{y}) \in \fml(T)$ 
such that both
\[
R_A[p(\bar{x},\bar{a})\union 
\phi(\bar{x},\bar{c})]\geq \alpha 
\quad
\text{ and }
\quad
R_A[p(\bar{x},\bar{a})\union 
\neg \phi(\bar{x},\bar{c})]\geq \alpha.
\] 
Since $A$ is $(D, \aleph_0)$-homogeneous, there is $\bar{d} \in A$ 
such that 
$\tp(\bar{c}\bar{a}, \empty)=\tp(\bar{d}\bar{b}, \empty)$. 
Therefore by induction hypothesis,
both 
\[
R_A[p(\bar{x},\bar{b})\union 
\phi(\bar{x},\bar{d})]\geq \alpha 
\quad
\text{ and }
\quad
R_A[p(\bar{x},\bar{b})\union 
\neg \phi(\bar{x},\bar{d})]\geq \alpha. 
\]
Second, for every $\bar{d} \in A$, there is 
$\bar{c} \in A$ such that
$\tp(\bar{c}\bar{a}, \empty)=\tp(\bar{d}\bar{b}, \empty)$. 
Thus, since $R_A[p(\bar{x},\bar{a})]\geq \alpha + 1$, 
there is $q(\bar{x}, \bar{y}) \in D$, such that
$R_A[p(\bar{x},\bar{a}) \cup q(\bar{x},\bar{c})] \geq \alpha$. 
Therefore, by induction hypothesis, 
$R_A[p(\bar{x},\bar{b}) \cup q(\bar{x},\bar{d})] \geq \alpha$.
This shows that $R_A[p(\bar{x},\bar{b})] \geq\alpha + 1$.
\end{itemize}

(4) 
Suppose $p \vdash q$. By definition of the
rank, we may choose
$q_0 \subseteq q$ over a finite set,
such that $R_A[q_0]=R_A[q]$.
Hence, since  $p \vdash q_0$, it is enough to show the
lemma when $q$ is over a finite set also.
Write $p=p(\bar{x},\bar{b}) \vdash q=q(\bar{x},\bar{a})$. 
We show by induction on
$\alpha$ that for every such pair of types over finite sets, 
we have
\[
R_A[p(\bar{x},\bar{b})] \geq \alpha
\quad
\text{ implies }
\quad 
R_A[q(\bar{x},\bar{b})] \geq \alpha.
\]
\begin{itemize}
\item 
For $\alpha=0$, this is true by definition.
\item
 For $\alpha$ a limit ordinal, this is true by induction.
\item 
Suppose $R_A[p(\bar{x},\bar{b})] \geq \alpha +1$.  
On the one hand, there is $\bar{c} \in A$ and 
$\phi(\bar{x},\bar{y}) \in \fml(T)$ such that both
\[
R_A[p(\bar{x},\bar{b})\union 
\phi(\bar{x},\bar{c})]\geq \alpha 
\quad
\text{ and }
\quad
R_A[p(\bar{x},\bar{b})\union 
\neg \phi(\bar{x},\bar{c})]\geq \alpha.
\] 
But 
\[
p(\bar{x},\bar{b})\union \phi(\bar{x},\bar{c}) \vdash
q(\bar{x},\bar{a})\union \phi(\bar{x},\bar{c})
\]
and similarly
\[
p(\bar{x},\bar{b})\union \neg \phi(\bar{x},\bar{c}) \vdash
q(\bar{x},\bar{a})\union \neg \phi(\bar{x},\bar{c}),
\] 
so by induction
hypothesis, both 
\[
R_A[q(\bar{x},\bar{a})\union 
\phi(\bar{x},\bar{c})]\geq \alpha 
\quad
\text{ and }
\quad
R_A[q(\bar{x},\bar{a})\union \neg 
\phi(\bar{x},\bar{c})]\geq \alpha.
\]
On the other hand, given any $\bar{c} \in A$, there is
 $r(\bar{x}, \bar{y}) \in D$, such that
$R_A[p(\bar{x},\bar{b}) \cup r(\bar{x},\bar{c})] \geq \alpha$.
But 
\[
p(\bar{x},\bar{b}) \cup r(\bar{x},\bar{c}) \vdash
q(\bar{x},\bar{a}) \cup r(\bar{x},\bar{c}),
\] 
so by induction hypothesis,
$R_A[q(\bar{x},\bar{a}) \cup r(\bar{x},\bar{c})] \geq \alpha$.
Hence $R_A[q(\bar{x},\bar{a})] \geq \alpha +1$.
\end{itemize}

(5)
First, choose $q(\bar{x}, \bar{a}) \subseteq p$, 
such that $R_A[q]=R_A[p]$ 
(this is possible by
definition of the rank). 
Similarly, since $f(q) \subseteq f(p)$, 
we could have chosen $q$ so that in addition $R_{f(A)}[f(q)]=R_{f(A)}[f(p)]$.
Now, by symmetry, it is enough to show that if
$R_A[q] \geq \alpha$ then $R_{f(A)}[f(q)]\geq \alpha$.
\begin{itemize}
\item 
For $\alpha = 0$ or $\alpha$ a limit ordinal, it is obvious by definition.
\item 
Suppose $\alpha = \beta + 1$. 
First,  there exists $\phi(\bar{x}, \bar{b})$
such that 
\[
R_A[q \cup \phi(\bar{x}, \bar{b})] \geq \beta
\quad
\text{and}
\quad
R_A[q \cup \neg \phi(\bar{x}, \bar{b})] \geq \beta.
\] 
Thus, by induction hypothesis,
we have  
\[
R_{f(A)}[f(q) \cup \phi(\bar{x}, \bar{f(b)})] \geq \beta
\quad 
\text{and}
\quad
R_{f(A)}[f(q) \cup \neg \phi(\bar{x}, \bar{f(b)})] \geq \beta.
\]
Second, notice that for every $\bar{b} \in f(A)$, 
there is $\bar{c} \in A$, such that
$f(\bar{c})=\bar{b}$. 
Since $R_A[q]\geq \beta + 1$, there exists $r(\bar{x}, \bar{y}) \in D$, 
such
that $R_A[q \cup r(\bar{x}, \bar{c})] \geq \beta$.
Hence, by induction hypothesis, 
$R_{f(A)}[f(q) \cup r(\bar{x}, \bar{b})] \geq \beta$.
This shows that 
$R_{f(A)}[f(q)] \geq \beta+1$. 
\end{itemize}

(6) 
This is immediate by definition of the rank.

(7)
By definition of the rank, let $B \in \dom(p)$ and 
$q \subseteq p \restriction B$
be such that $R_A[q]=R_A[p]$.
Now, clearly $q \subseteq p \restriction B \subseteq p$, 
so $R_A[q] \geq R_A[p \restriction B] \geq R_A[p]$ by Lemma \ref{subset}.
So $R_A[p\restriction B]=R[p]$.

(8)
Suppose there is $\alpha_0$ such that $R_A[p] \not = \alpha_0$ for every $p$.
We prove by induction on $\alpha \geq \alpha_0$, 
that for no type $p$ do we have $R_A[p]= \alpha$. 
\begin{itemize}
\item 
For $\alpha = \alpha_0$, this is the definition of $\alpha_0$.
\item 
Now suppose that there is $p$ such that $R_A[p]=\alpha +1$.
By \ref{finite}, we may assume that $p$ is over a finite set.
Then there is $\bar{c} \in A$ and $\phi(\bar{x},\bar{y}) \in \fml(T)$
 such that both 
\[
R_A[p \union \phi(\bar{x},\bar{c})]\geq \alpha 
\quad
\text{ and }
\quad
R_A[p \union \neg \phi(\bar{x},\bar{c})]\geq \alpha.
\] 
But by induction hypothesis, neither
can be equal to $\alpha$, so we must have both 
\[
R_A[p \union \phi(\bar{x},\bar{c})]\geq \alpha+1
\quad
\text{ and }
\quad
R_A[p \union \neg \phi(\bar{x},\bar{c})]\geq \alpha+1.
\] 
Similarly, given any $\bar{c} \in A$, 
there is $q(\bar{x}, \bar{y}) \in D$, 
such that $R_A[p \cup q(\bar{x},\bar{c})] \geq \alpha$. 
But, by induction hypothesis, we cannot
have $R_A[p \cup q(\bar{x},\bar{c})] = \alpha$, so 
$R_A[p \cup q(\bar{x},\bar{c})] \geq \alpha+1$. 
But this shows that $R_A[p] \geq \alpha + 2$, a
contradiction.
\item 
Suppose $\alpha > \alpha_0$ is a limit ordinal. 
Then $\alpha \geq \alpha_0 + 1$, so as in the 
previous case, there is $\bar{c} \in A$ and 
$\phi(\bar{x},\bar{y}) \in \fml(T)$
such that both 
\[
R_A[p \union \phi(\bar{x},\bar{c})]\geq \alpha_0
\quad
\text{ and }
\quad
R_A[p \union \neg \phi(\bar{x},\bar{c})]\geq \alpha_0.
\] 
But by induction hypothesis, for no
$\beta$ such that $\alpha > \beta  \geq \alpha_0$ can we have 
$R_A[p \union \phi(\bar{x},\bar{c})]=\beta$ or
$R_A[p \union \neg \phi(\bar{x},\bar{c})]=\beta$, 
so necessarily since $\alpha$ is a limit
ordinal, we have  
\[
R_A[p \union \phi(\bar{x},\bar{c})]\geq\alpha 
\quad
\text{ and }
\quad
R_A[p \union \neg \phi(\bar{x},\bar{c})]\geq\alpha.
\]
Similarly, for any $\bar{c} \in A$, 
there is $q(\bar{x}, \bar{y}) \in D$, 
such that $R_A[p \cup q(\bar{x},\bar{c})] \geq \alpha_0$ 
and hence by induction hypothesis
$R_A[p \cup q(\bar{x},\bar{c})] > \beta$ 
for any $\alpha_0 \leq \beta < \alpha$ 
so since $\alpha$ is a limit ordinal, 
we have $R_A[p \cup q(\bar{x},\bar{c})] \geq \alpha$. 
But this shows that $R_A[p] \geq \alpha + 1$, a
contradiction.
\end{itemize}

(9) 
By the previous lemma, it is enough to find $\alpha_0 <(|A| +2^{|T|})^+$,
(respectively $< (2^{|T|})^+$ if $A$ is a $(D, \aleph_0)$-homogeneous model) 
such that 
\begin{equation} \tag{*}
R_A[p] \not = \alpha_0 
\qquad 
\text{for every  type over $A$}.
\end{equation} 
We do this by counting the number of possible values for the rank.
By \ref{finite} it is enough to count 
the values achieved by types over finite subsets of $A$. 
But there are at most $|A|^{<\aleph_0}\leq |A| + \aleph_0$ finite subsets of 
$A$, and given any finite subset, there are only $2^{|T|}$ distinct types 
over it. 
Hence there are at most $|A| +2^{|T|}$ 
many different ranks, and so by the pigeonhole principle
(*) holds for some $\alpha_0 < ( |A| + 2^{|T|})^+$.

When $A$ is a $(D,\aleph_0)$-homogeneous model, 
the bound can be further reduced
by a use of \ref{2}, since only the type of
each of those finite subset of $A$ is relevant.
\epf

The next lemma shows that the rank is especially well-behaved
when the parameter $A$ is the universe of a $(D, \aleph_0)$-homogeneous model. 
This is used in particular to study 
$(D,\aleph_0)$-homogeneous models in the last two sections.

\begin{lemma} \label{homo}

\begin{enumerate} 
\item 
If $p$ is over a subset of a $(D,\aleph_0)$-homogeneous model $M$, 
then $R_M[p]=R[p].$

\item 
If $p$ is over $M_1 \intersect M_2$, with $M_l$ $(D,\aleph_0)$-homogeneous, 
for $l=1,2$,
we have $R_{M_1}[p]=R_{M_2}[p].$

\item 
If $q(\bar{x},\bar{a}_l)$ are sets of formulas, with $a_l \in M_l$ for
$l=1,2$ satisfying $\tp(\bar{a}_1, \empty)=\tp(\bar{a}_2,\empty)$, then
$R_{M_1}[q(\bar{x},\bar{a}_1)] =R_{M_2}[q(\bar{x},\bar{a}_2)].$

\end{enumerate}

\end{lemma}

\bpf(1)  
First, by Finite Character, we may assume that $p$ is over a finite set.
Now we show by induction on $\alpha$ that 
\[
R_M[p] \geq \alpha
\qquad
\text{implies}
\qquad
R[p] \geq \alpha.
\]
When $\alpha = 0$ or $\alpha$ is a limit, it is clear.
Suppose $R_M[p] \geq \alpha+1$. 
Then there is $\bar{b} \in M$ and 
$\phi(\bar{x},\bar{y})$ such that both
\[
R_M[p \cup \phi(\bar{x},\bar{b})] \geq \alpha
\quad
\text{ and }
\quad
R_M[p \cup \neg \phi(\bar{x},\bar{b})] \geq \alpha.
\]
By induction hypothesis, we have 
\[
R[p \cup \phi(\bar{x},\bar{b})] \geq \alpha 
\quad
\text{ and }
\quad
R[p \cup \neg \phi(\bar{x},\bar{b})] \geq \alpha.
\]
Further, if $\bar{b} \in \mathfrak{C}$, 
choose $\bar{b}' \in M$, such that
$\tp(\bar{b},\bar{a})=\tp(\bar{b}',\bar{a})$. 
Since $R_M[p]\geq \alpha+1$, there is 
$q(\bar{x},\bar{y}) \in D$ such that
$R_M[p \cup q(\bar{x},\bar{b}')]\geq \alpha.$
Thus, since $\mathfrak{C}$ is $(D,\aleph_0)$-homogeneous,  
by induction hypothesis
we have 
$R[p \cup q(\bar{x},\bar{b}')]\geq \alpha$, and so
by Lemma \ref{rankproperties} \ref{2} 
$R[p \cup q(\bar{x},\bar{b})]\geq \alpha$. 
Hence $R[p]\geq \alpha+1$.

For the converse, similarly by induction on $\alpha$ we show that 
\[
R[p] \geq \alpha
\qquad
\text{implies}
\qquad
R_M[p] \geq \alpha.
\]
Again, for $\alpha=0$ or $\alpha$ a limit, it is easy.
Suppose $R[p] \geq \alpha+1$. Then there is $\bar{b} \in \mathfrak{C}$ and 
$\phi(\bar{x},\bar{y})$ such that both
\[
R[p \cup \phi(\bar{x},\bar{b})] \geq \alpha
\quad
\text{ and }
\quad
R[p \cup \neg \phi(\bar{x},\bar{b})] \geq \alpha.
\]
Since $M$ is $(D,\aleph_0)$-homogeneous, there exists
 $\bar{b}' \in M$, such that
$\tp(\bar{b},\bar{a})=\tp(\bar{b}',\bar{a})$. 
By Lemma \ref{rankproperties} \ref{2}, we have 
\[
R[p \cup \phi(\bar{x},\bar{b}')] \geq \alpha
\quad
\text{ and }
\quad
R[p \cup \neg \phi(\bar{x},\bar{b}')] \geq \alpha.
\]
Hence, by induction hypothesis, we have (since $\bar{b}'\in M$)
\[
R_M[p \cup \phi(\bar{x},\bar{b}')] \geq \alpha
\quad
\text{ and }
\quad
R_M[p \cup \neg \phi(\bar{x},\bar{b}')] \geq \alpha.
\]
Also, for any $\bar{b} \in M$, since $\bar{b} \in \mathfrak{C}$
there is $q(\bar{x},\bar{y}) \in D$ such that
$R[p \cup q(\bar{x},\bar{b})]\geq \alpha.$
By induction hypothesis, we have 
$R_M[p \cup q(\bar{x},\bar{b})]\geq \alpha,$
which finishes to show that $R_M[p]\geq \alpha+1$ and completes the proof.

(2)
By (1) applied twice, $R_{M_1}[p]=R[p]=R_{M_2}[p]$.

(3)
Since $R_{M_1}[q(\bar{x},\bar{a}_1)] 
=R[q(\bar{x},\bar{a}_1)]=
R[q(\bar{x},\bar{a}_2)]=
R_{M_2}[q(\bar{x},\bar{a}_2)]$.
\epf

We now show that the rank is bounded when $D$ is $\aleph_0$-stable.

\begin{theorem} \label{existence}
If $D$ is stable in $\lambda$ for some $\aleph_0 \leq \lambda < 2^{\aleph_0}$
then  
$R_A[p] < \infinity$ for every type $p$ and every subset $A$ of $\mathfrak{C}$.
\end{theorem}
\bpf 
We prove the contrapositive.
Suppose there is a subset $A$ of $\mathfrak{C}$ and a type 
$p$ over $A$ 
such that $R_A[p]=\infinity$.
We construct sets $A_{\eta} \subseteq A$ and 
types $p_{\eta}$, for $\eta \in \sq{<\omega}2$, 
such that:
\begin{enumerate}

\item 
$p_{\eta} \in S_D(A_{\eta})$;

\item 
$p_{\eta} \subseteq p_{\nu}$ 
when $\eta < \nu$;

\item 
$A_{\eta}$ is finite;

\item 
$p_{\eta \conc 0}$ and $p_{\eta \conc 1}$ are contradictory;

\item 
$R_A[p_{\eta}]= \infinity$;

\end{enumerate}

This is possible: 
Let $\mu = (2^{|T|})^+$ if $A$ is a $(D,\aleph_0)$-homogeneous model,
and $\mu=(|A|+2^{|T|})^+$ otherwise. The construction is by 
induction on $n=\ell(\eta)$.
\begin{itemize}
\item
For $n=0$, by Finite Character we choose first $\bar{b} \in A$, 
such that $R_A[p]=R_A[p \restriction \bar{b}]=\infinity$.
Since $R_A[p \restriction \bar{b}]=\infinity$, in particular 
$R_A[p \restriction \bar{b}] \geq \mu +1$ so
there exists $q(\bar{x},\bar{y})\in D$, such that 
$R_A[(p\restriction \bar{b})\cup q(\bar{x},\bar{b})]\geq \mu.$ 
But then $p\restriction \bar{b}\subseteq q(\bar{x},\bar{b})$, 
$q(\bar{x},\bar{b}) \in S_D(\bar{b})$ and 
$R_A[q(\bar{x},\bar{b})]\geq \mu$, so
$R_A[q(\bar{x},\bar{b})] =\infinity$
 by Lemma \ref{rankproperties} \ref{bound}.
Therefore, we let $A_{<>}=\bar{b}$ and 
$p_{<>}=q(\bar{x},\bar{b})$ and the conditions are satisfied.

\item
Assume $n \geq 0$ and that we have
constructed $p_{\eta} \in S_D(A_{\eta})$ with $\ell(\eta)=n$. 
Since $R_A[p_{\eta}]= \infinity$, in particular 
$R_A[p_{\eta}] \geq (\mu +1) +1.$
Hence, 
there is $\bar{a}_{\eta} \in A$
and $\phi(\bar{x}, \bar{y})$ such that 
\begin{equation} \tag{*}
R_A[p_{\eta} \cup \phi(\bar{x}, \bar{a}_{\eta})]\geq \mu+1 
\quad
\text{ and }
\quad
R_A[p_{\eta} \cup \neg \phi(\bar{x}, \bar{a}_{\eta})]\geq \mu +1.
\end{equation}
Let $A_{\eta \conc 0}=A_{\eta \conc 1}=
A_{\eta} \cup \bar{a}_{\eta} \subseteq A$. 
Both $A_{\eta \conc 0}$ and $A_{\eta \conc 1}$ are finite, 
so (*) and the definition of the rank imply 
that 
there are $q_l(\bar{x},\bar{y})\in D$ for $l=0,1$, such that 
\[
R_A[p_{\eta} \cup \phi(\bar{x}, \bar{a}_{\eta}) \cup 
q_0(\bar{x},A_{\eta \conc 0})] \geq \mu
\]
 and
\[
R_A[p_{\eta} \cup \neg \phi(\bar{x}, \bar{a}_{\eta}) 
\cup q_1(\bar{x},A_{\eta \conc 1})] \geq \mu.
\]
Define $p_{\eta \conc 0}:=p_{\eta} \cup \phi(\bar{x}, \bar{a}_{\eta}) 
\cup q_0(\bar{x},A_{\eta \conc 0})$ and
$p_{\eta \conc 1}:=p_{\eta} \cup \neg \phi(\bar{x}, \bar{a}_{\eta}) 
\cup q_1(\bar{x},A_{\eta \conc 1})$.
Then $p_{\eta \conc l} \in S_D(A_{\eta \conc l})$ since 
$q_l(\bar{x},A_{\eta \conc l}) \in S_D(A_{\eta \conc l})$ and 
$A_{\eta \conc l}$ is finite for $l=0,1$. 
Moreover,  
$p_{\eta \conc 0}$ and $p_{\eta \conc 1}$ are contradictory by construction.
Finally
$R_A[p_{\eta \conc l}]= \infinity$, 
since $R_A[p_{\eta \conc l}]\geq \mu$. 
Hence all the requirements are met.

\end{itemize}
This is enough: For each $\eta \in \sq\omega2$, define 
$A_{\eta}:= \Union_{n \in \omega} A_{\eta \restriction n}$ and
$p_{\eta}:= \Union_{n \in \omega} p_{\eta \restriction n}$.
We claim that $p_{\eta} \in S_D(A_{\eta})$. 
Certainly
$p_{\eta} \in S(A_{\eta})$, so we only need to 
show that if $\bar{c} \models p_{\eta}$, then
$A_{\eta} \cup \bar{c}$ is a $D$-set ($\bar{c}$ is not assumed
to be in $\mathfrak{C}$).
It is enough to show that
$\tp(\bar{c}\bar{d}, \empty) \in D$ for every finite
$\bar{d} \in A_{\eta}$.
But, if $\bar{d} \in A_\eta$, then there is
$n \in \omega$ such that $\bar{d} \in A_{\eta \restriction n}$. 
Since 
$\bar{c} \models p_{\eta \restriction n}$ and 
$p_{\eta \restriction n} \in S_D(A_{\eta \restriction n})$, then
$\bar{c} \cup A_{\eta \restriction n}$ is a $D$-set, and
therefore
$\tp(\bar{c}\bar{d}, \empty) \in D$, which is what we wanted.
Now that we have established that
$p_\eta \in S_D(A_\eta)$,  since $\mathfrak{C}$ is $(D, \chi)$-homogeneous, 
there is
$\bar{c}_{\eta} \in \mathfrak{C}$ such
that $\bar{c}_{\eta} \models p_{\eta}$. 
Now let $C= \Union_{ \eta \in \sq{<\omega} 2} A_{\eta}$. 
Then $|C|=\aleph_0$ and if $\eta \not = \nu \in \sq\omega2$, 
then $\tp(\bar{c}_\eta,C)\not = \tp(\bar{c}_\nu,C)$,
since $p_{\eta}$ and 
$p_{\nu}$ are contradictory.
Therefore  $|S_D(C)|\geq 2^{\al{0}}$, which
shows that $D$ is not stable in $\lambda$ for any 
$\aleph_0 \leq \lambda < 2^{\aleph_0}$.
\epf

\begin{remark} 
Recall that in \cite{Sh:1}, $D$ is stable in $\lambda$ if and only if 
there is a $(D,\lambda^+)$-homogeneous model and 
$|S_D(A)|\leq \lambda$ for all $D$ sets $A$ of cardinality at most $\lambda$ 
(this is Definition 2.1 of \cite{Sh:1}). 
The proof of the previous theorem shows
that if $D$ is stable in $\lambda$ for some 
$\aleph_0 \leq \lambda < 2^{\aleph_0}$ in the sense of \cite{Sh:1}
then $R_A[p] < \infinity$
for all $D$-set $A$ and $D$-type $p$. 
In other words, we do not really need $\mathfrak{C}$ for this proof.
\end{remark}

By analogy with the first-order case (see \cite{Sh:a} definition 3.1), 
we introduce the following definition:

\begin{definition} We say that $D$ is \emph{totally transcendental} if 
$R_A[p] < \infinity$ for every subset $A$ of $\mathfrak{C}$
and every type $p$ over $A$.
\end{definition}

For the rest of the paper, we will make
the following hypothesis. We will occasionally repeat that $D$
is totally transcendental for emphasis.
\begin{hypothesis} $D$ is totally transcendental.
\end{hypothesis}

In what follows, we shall show that when $D$ is totally transcendental, 
the rank affords a well-behaved dependence relation on the
subsets of $\mathfrak{C}$.
We first focus on a special kind of types.

\begin{definition} \label{stationarydef} 
A type $p$ is called \emph{stationary} if for every $B$ containing  $\dom(p)$ 
there is a unique type $p_B \in S_D(B)$, 
such that $p_B$ extends $p$ and $R[p]=R[p_B]$.
\end{definition}

Note that since our rank is not an extension of Morley's rank, 
one does not necessarily
get the usual stationary types when the class is first-order. 
The argument in the next lemma is a generalization of 
Theorem 1.4.(1)(b) in \cite{Sh:871}.
Recall that $p \in S_D(A)$ \emph{splits over $B \subseteq A$}
if there exists
$\phi(\bar{x}, \bar{y})$ and $\bar{a}, \bar{c} \in A$ with
$\tp(\bar{a},B)=\tp(\bar{c},B)$, 
such that  $\phi(\bar{x}, \bar{a}) \in p$ and 
$\neg \phi(\bar{x}, \bar{c}) \in p$.

\begin{lemma} \label{exten}
Suppose there is  $\bar{d} \in \mathfrak{C}$ realizing 
$p(\bar{x}, \bar{b})$ and a  $(D,\aleph_0)$-homogeneous model
$M$ 
such that
\begin{equation} \tag{*}
R[\tp(\bar{d},M)]=R[p(\bar{x},\bar{b})]=\alpha.
\end{equation} 
Then, for any $A \subseteq \mathfrak{C}$ containing $\bar{b}$ 
there is a unique $p_A \in S_D(A)$ extending
$p(\bar{x},\bar{b})$, such that 
\[
R[p_A]=R[p(\bar{x},\bar{b})]=\alpha.
\]
Moreover, $p_A$ does not split over $\bar{b}$.
\end{lemma}

\bpf We first prove uniqueness.
Suppose two different types $p_A$ and $q_A \in S_D(A)$ extend 
$p(\bar{x},\bar{b})$ and 
\[
R[p_A]=R[p(\bar{x},\bar{b})]=R[q_A]=\alpha.
\]
Then there is $\phi(\bar{x},\bar{c}) \in p_A$ such that 
$\neg \phi(\bar{x},\bar{c}) \in q_A$.
Thus, by Monotonicity,
\[ 
R[p(\bar{x},\bar{b}) 
\cup \phi(\bar{x},\bar{c})]
\geq
R_A[p] = \alpha
\quad
\text{ and }
\quad
R[p(\bar{x},\bar{b}) 
\cup \neg \phi(\bar{x},\bar{c})]
\geq 
R_A[p]= \alpha.
\]
Further, for every $\bar{c} \in \mathfrak{C}$, 
there is $\bar{c}' \in M$ such that 
$\tp(\bar{c},\bar{b})=\tp(\bar{c}',\bar{b})$ since $M$ is 
$(D,\aleph_0)$-homogeneous.
Now write $q(\bar{x}, \bar{c}')=\tp(\bar{d}, \bar{c}')$, and notice that
\[
R[p(\bar{x},\bar{b}) \cup q(\bar{x}, \bar{c}')]
\geq
R[\tp(\bar{d}, \bar{b} \cup \bar{c}')]                        
\geq   
R[\tp(\bar{d},M)]=\alpha. 
\] 
But $q(\bar{x}, \bar{y}) \in D$ 
by definition and so by Lemma \ref{rankproperties} (2)
$R[p(\bar{x},\bar{b}) \cup q(\bar{x}, \bar{c})] \geq \alpha$ 
since 
$\tp(\bar{c}\bar{b},\empty)=\tp(\bar{c}'\bar{b},\empty)$. 
But this shows that 
$R[p(\bar{x},\bar{b})] \geq \alpha+1$, which contradicts (*).

We now argue that $p_A$ does not 
split over $\bar{b}$. Suppose it does, and choose
a formula $\phi(\bar{x},\bar{y}) \in \fml(T)$ 
and sequences $\bar{c}_0, \bar{c}_1 \in A$
with $\tp(\bar{c}_0,\bar{b})=\tp(\bar{c}_1,\bar{b})$ 
such that 
$\phi(\bar{x},\bar{c}_0) $ and 
$\neg \phi(\bar{x},\bar{c}_1)$ both belong to $p_A$.
Then by Monotonicity,
\[
R[p(\bar{x},\bar{b}) \cup \phi(\bar{x},\bar{c}_0)]
\geq
R_A[p]=\alpha
\quad
\text{ and }
\quad
R[p(\bar{x},\bar{b}) \cup \neg \phi(\bar{x},\bar{c}_1)]
\geq
R_A[p]=\alpha.
\]
But $\tp(\bar{c}_0,\bar{b})=\tp(\bar{c}_1,\bar{b})$ 
so by Lemma \ref{rankproperties}(3) we have
\[
R[p(\bar{x},\bar{b}) \cup \phi(\bar{x},\bar{c}_1)] \geq \alpha.
\]
 An argument similar to the uniqueness argument in the first paragraph
finishes to show that $R[p(\bar{x},\bar{b})] \geq \alpha+1$, which is
again
a contradiction to (*).

For the existence, let $p_A$ be the following
set of formulas with parameters
in $A$:
\[
\{\, \phi(\bar{x},\bar{c}) \mid \,  
\text{There exists }
\bar{c}' \in M
\text{ such that }
\tp(\bar{c},\bar{b})=\tp(\bar{c}',\bar{b})
\text{ and }
\models \phi[\bar{d},\bar{c}'] \,\}.
\]
By the non-splitting part, using the fact that 
$M$ is $(D,\aleph_0)$-homogeneous, 
we have that $\tp(\bar{d},M)$ does not split over $\bar{b}$.
Hence $p_A \in S_D(A)$ and 
does not split over $\bar{b}$.
We show that this implies that 
$R[p_A]=R[\tp(\bar{d},M)]=\alpha.$ 
Otherwise, since $p_A$ extends $p(\bar{x},\bar{b})$, by Monotonicity
we must have $R[p_A] \leq \alpha$,
and therefore $R[p_A]< \alpha$.
Let us choose $\bar{b}' \in A$ such that $\bar{b} \subseteq \bar{b}'$ and 
$R[p_A]=R[p_A \restriction \bar{b}']$. 
For convenience, we write $q(\bar{x},\bar{b}') := p_A \restriction \bar{b}'$,
and so $R[q(\bar{x},\bar{b}')] < \alpha$.
Now since $M$ is $(D,\aleph_0)$-homogeneous, we can choose $\bar{b}'' \in M$ 
such that 
$\tp(\bar{b}'', \bar{b})=\tp(\bar{b}', \bar{b})$. 
Hence 
\begin{equation} \tag{**}
R[q(\bar{x}, \bar{b}')]= R[q(\bar{x},\bar{b}')] < \alpha.
\end{equation}
But by definition of $p_A$, we must have 
$q(\bar{x},\bar{b}') \subseteq \tp(\bar{d},M)$, 
so by Monotonicity we have
$R[q(\bar{x},\bar{b}')] \geq R[ \tp(\bar{d},M)] = \alpha$,
which contradicts (**).
\epf

\begin{corollary}\label{stationary}
The following conditions are equivalent:
\begin{enumerate}
\item
$p \in S_D(A)$ is stationary. 
\item
There is a 
$(D,\aleph_0)$-homogeneous model $M$ containing $A$
and $\bar{d} \in \mathfrak{C}$ realizing $p$
such that $R[\tp(\bar{d}, M)]=R[p]$.
\end{enumerate}
\end{corollary}

\begin{definition} A stationary type $p \in S_D(A)$ is \emph{based on $B$} if
$R[p]=R[p \restriction B]$.
\end{definition}

\begin{remark}
\begin{enumerate}
\item
If $p$ is stationary, 
there is a finite $B \subseteq \dom(p)$ such that $p$ is based on $B$.
\item 
If $p$ is based on $B$, then $p \restriction B$ is also stationary
and $p$ is the only extension of $p\restriction B$ such that
$R[p]=R[p\restriction B]$. 
\item 
If $p$ is stationary and $\dom(p) \subseteq A \subseteq B$,
then $p_A = p_B \restriction A$.
\item 
Suppose $\tp(\bar{a}, \empty)=\tp(\bar{a}', \empty)$. 
Then $p(\bar{x},\bar{a}')$ is stationary if and only if 
$p(\bar{x},\bar{a})$ is stationary.
(Use an automorphism of $\mathfrak{C}$ sending $\bar{a}$ to $\bar{a}'$.) 
\end{enumerate}
\end{remark}

Stationary types allow us to prove a converse of Theorem \ref{existence}.

\begin{theorem} \label{stable} 
If $D$ is totally transcendental then $D$ is
stable in every $\lambda \geq |D|+|T|$. In particular $\kappa(D)=\aleph_0$.
\end{theorem}

\bpf 
Let $\lambda \geq |D|+|T|$, and let $A$ be a subset of $\mathfrak{C}$ 
of cardinality at most $\lambda$.
Since $\lambda \geq |D|+|T|$, by using a countable, increasing chain of models
we can find a $(D, \aleph_0)$-homogeneous model $M$ containing
$A$ of cardinality 
$\lambda$.
Since $|S_D(A)| \leq |S_D(M)|$, it is enough to show that $|S_D(M)| \leq \lambda$.
Suppose that $|S_D(M)| \geq \lambda^+$. 
Since $M$ is $(D, \aleph_0)$-homogeneous, each $p \in S_D(M)$ is stationary.
Hence, for each $p \in S_D(M)$, 
 we can choose a finite $B_p \subseteq M$ 
such that $p$ is based on $B_p$.
Since there are only $\lambda$ many
finite subsets of $M$, by the pigeonhole principle
 there is a fixed finite subset $B$ of $M$ 
such that $\lambda^+$ many types $p \in S_D(M)$ are based on $B$.
Since $\lambda^+ > |S_D(B)|=|D|$, another application of
the pigeonhole principle shows that there a single stationary type
$q \in S_D(B)$ with  $\lambda^+$ many extensions in $S_D(M)$
of the same rank. This contradicts the stationarity of $q$.
Hence $D$ is stable in $\lambda$.

For the last sentence, let $\lambda = \beth_\omega (|D|+|T|)$. 
By Zermelo-K\"{o}nig,
$\lambda^{\aleph_0} > \lambda$, hence by Theorem \ref{stabs} 
$\kappa(D)=\aleph_0$.
\epf

The following results show that stationary types behave nicely. Not only do
they have the uniqueness and the extension properties, 
but they can be represented by averages.
Surprisingly, it turns out that every type is reasonably 
close to a stationary type (this is made precise in Lemma \ref{nice}).

\begin{definition} Let $p \in S_D(A)$ be stationary
and let $\alpha$ be an infinite ordinal.
The sequence
$I = \{\, c_i \mid i < \alpha \,\}$ is called a 
\emph{Morley sequence based on $p$} if
for each $i < \alpha$ we have
$c_i$ realizes $p_{A_i}$, 
where 
$A_i = A \cup \{ c_j \mid j<i \}$.
\end{definition}

\begin{lemma} \label{morley} Let $p\in S_D(A)$ be stationary.
If $I$ is a Morley sequence based on $p$, then $I$ is indiscernible over $A$.
\end{lemma}

\bpf By stationarity $p_{A_i} \subseteq p_{A_j}$ when $i<j$, and
by the previous lemma each $p_{A_i}$ does not split over $A$.
Hence, a standard result (see for example \cite{Sh:a} Lemma I.2.5) 
implies 
that $I$ is an indiscernible sequence over $A$.
\epf

\begin{definition} ($\kappa(D)=\aleph_0$)
For $I$ an infinite set of indiscernibles and $A$ a set 
(with $I \cup A \subseteq \mathfrak{C}$),
recall that
\[
\av(I,A)=\{\, \phi(\bar{x}, \bar{a}) \mid
\bar{a} \in A, \phi(\bar{x},\bar{y}) \in L(T)
 \text{ and } 
|\phi(I,\bar{a})|\geq \aleph_0 \,\}.
\]
\end{definition}

\begin{lemma} \label{averagemorley} Suppose 
$p \in S_D(A)$ is stationary and
 $I$ is a Morley sequence based on $p$.
Then for any $B$ containing $A$
 we have that $p_B=\av(I,B)$.
\end{lemma}

\bpf 
Let $B \subseteq \mathfrak{C}$ and write $I = \{ c_i \mid i < \alpha \}$.
Choose $c_i \in \mathfrak{C}$ for 
$\alpha \leq i < \alpha+\omega$ realizing $p_{B_i}$,
where $B_i=B \cup \Union \{a_j \mid j<i \}$.
Since $\av(I,B) \in S_D(B)$ extends $p$, it is enough to show
that $R[\av(I,B)]=R[p]$. 
Suppose $R[\av(I,B)]\not = R[p]$. 
Then, by Monotonicity, we must have $R[\av(I,B)] < R[p]$. 
We can find a finite $C \subseteq B$ such that $p$ is based
on $C$ and by Finite Character, we may assume in addition that 
\begin{equation} \tag{*}
R[\av(I,B)]=R[\av(I,C)] < R[p].
\end{equation}
But, since $C$ is finite and $\kappa(D)=\aleph_0$, by Lemma \ref{kappa} 
there is  $c_i \in I$ for $\alpha \leq i <\alpha+\omega$
realizing $\av(I,C)$, and since $C \subseteq B$, 
we must have $\tp(c_i, C)=\av(I,C)=p_C$
(since $c_i$ realizes $p_{B_i}$).
But then, by choice of $C$ we have $R[\av(I,C)]=R[p_C]=R[p]$
which contradicts (*).  
\epf

\begin{lemma} 
Let $I$ be an infinite indiscernible set, $A$ be finite 
and $p=\av(I,A)$ be stationary.
Then for any $C \supseteq A$ we have $p_C=\av(I,C)$.
\end{lemma}

\bpf Write $I=\{c_i \mid i<\alpha\}$, for $\alpha \geq \omega$ and 
let $C$ be given.
Choose $c_i \in \mathfrak{C}$ for $\alpha \leq i < \alpha+\omega$ 
realizing $p_{C_i}$, where $C_i=C\cup \Union \{c_j \mid j<i\}$.
Let $I'=\{c_i \mid i< \alpha+\omega \}$ and notice 
that necessarily $\av(I,B)=\av(I',B)$ for any $B$.
Suppose $p_C \not = \av(I,C)$, then since 
$\av(I,A) \subseteq \av(I,C)$, we must have
$R[\av(I,C)]<R[p]$, so $R[\av(I',C)]<R[p_C]$. 
Choose $C'$ finite, with $A \subseteq C' \subseteq C$, 
such that $R[\av(I',C)]=R[\av(I',C')]$.
Now there is $J \subseteq I'$ finite such that $I'-J$ 
is indiscernible over $C'$.
Choose $c_i \in I'-J$ with $i>\alpha$. 
Then $c_i$ realizes $\av(I',C')$, so 
$\av(I',C')=\tp(c_i,C')\subseteq p_{C_i}$ by choice of $c_i$.
But then 
\[
R[\av(I',C')]\geq R[p_{C_i}]=R[p]>R[\av(I,C)]=R[\av(I',C')], 
\]
a contradiction. 
\epf

It is natural at this point to introduce the forking symbol, 
by analogy with the
first-order case (see for example
\cite{Bl} or \cite{Ma}). 
We do not claim that the two notions coincide even when both are defined.

\begin{definition} \label{fork}
Suppose $A$, $B$, $C \subseteq \mathfrak{C}$, with $B \subseteq A$.
We say that 
\[
A \nonfork_B C 
\quad
\text{ if }
\quad
R[\tp(\bar{a}, B)]=R[\tp(\bar{a},B \cup C)],
\qquad \text{for every $\bar{a} \in A$}.
\]
\end{definition}

As in many other contexts, the symmetry property can be obtained
from the failure
of the order property.

\begin{theorem} [Symmetry]\label{symmetry} 
If $\tp(\bar{a},B)$ and $\tp(\bar{c},B)$ are stationary,
then 
\[
\bar{a} \nonfork_B \bar{c} 
\qquad 
\text{if and only if}
\qquad
\bar{c} \nonfork_B \bar{a}.
\]
\end{theorem}

\bpf First, $D$ is stable  by Theorem \ref{stable},
and therefore does not have the $\infinity$-order property 
by Theorem \ref{unstable}. Suppose, for a contradiction, that
\[
R[\tp(\bar{c}, B \cup \bar{a})] < R[\tp(\bar{c}, B)]
\quad
\text{ and }
\quad
R[\tp(\bar{a}, B \cup \bar{c})]=R[\tp(\bar{a}, B)].
\] 
Let $\lambda = \beth_{(2^{|T|})^+}$ and let $\mu =(2^\lambda)^+$. 
We use Theorem \ref{order} to show that $D$ has the $\infinity$-order property,
by constructing an order of length $\lambda$.
Choose $p(\bar{x},\bar{y}, \bar{b}) \in S_D(\bar{b})$ with 
$\bar{b} \in B$, such that
\[
R[\tp(\bar{a}, B \cup \bar{c})]=
R[p(\bar{x},\bar{c}, \bar{b})]=
R[\tp(\bar{a}, B)]
\]
and
\[
R[\tp(\bar{c}, B\cup \bar{a})]=
R[p(\bar{c},\bar{y}, \bar{b})]<
R[\tp(\bar{c}, B)].
\]
Let $\bar{a}_{\alpha}, \bar{c}_{\alpha} \in \mathfrak{C}$ 
for $\alpha < \mu$ 
and 
$B_{\alpha}= \Union \{ \bar{a}_\beta, \bar{c}_\beta \mid
 \beta < \alpha \}$ be such that:

\begin{enumerate}
\item 
$B_0 = B$;
\item 
$\bar{a}_{\alpha}$ realizes $\tp(\bar{a},B)$
and 
$R[\tp(\bar{a}_{\alpha}, B_{\alpha})]=R[\tp(\bar{a},B)]$;

\item 
$\bar{c}_{\alpha}$ realizes $\tp(\bar{c},B)$ and
$R[\tp(\bar{c}_{\alpha}, B_{\alpha}\cup \bar{a}_\alpha)]=
R[\tp(\bar{c},B)]$.

\end{enumerate}

This is achieved by induction on $\alpha < \mu$.
Let $B_0:=B$, $\bar{a}_0:=\bar{a}$
and $\bar{c}_0:=\bar{c}$. At stage $\alpha$, we
let first 
$B_{\alpha}:= \Union \{ \bar{a}_\beta, \bar{c}_\beta \mid \beta < \alpha  \}$
which is well-defined by induction hypothesis. 
We then satisfy in this order
(2) by stationarity of $\tp(\bar{a}, B)$,
and (3) by stationarity of $\tp(\bar{c}, B)$. 

This is enough: First, notice that $\bar{c}_{\alpha}$ does not realize 
$p(\bar{a},\bar{y}, \bar{b})$, otherwise 
\[
R[\tp(\bar{c}_{\alpha}, B_{\alpha}\cup \bar{a}_\alpha)] 
\leq R[p(\bar{a},\bar{y}, \bar{b})]
<R[\tp(\bar{c}, B)],
\] 
contrary to the choice of $\bar{c}_{\alpha}$. 
Similarly, since $\tp(\bar{a}_{\alpha},B)=\tp(\bar{a},B)$ and
$\bar{b} \in B$, then 
\[
R[p(\bar{a}_{\beta},\bar{y}, \bar{b})]<R[\tp(\bar{c}, B)],
\] 
so $\bar{c}_{\alpha}$ does not realize 
$p(\bar{a}_{\beta},\bar{y}, \bar{b})$ when 
$\alpha \geq \beta$.

Now suppose $\alpha < \beta$. 
Then $\bar{a}_{\beta}$ realizes 
$p(\bar{x},\bar{c},\bar{b})$ since by stationarity,
 we must have 
$\tp(\bar{a}_{\beta}, A \cup \bar{c})= 
\tp(\bar{a}, B\cup \bar{c})$.
Further, since $\tp(\bar{a}_{\alpha}, B_{\alpha})$ does not split over $B$ and 
$\tp(\bar{c}_{\alpha},B)=\tp(\bar{c},B)$ we must have 
$p(\bar{x},\bar{c}_{\alpha},\bar{b}) 
\subseteq \tp(\bar{a}_{\alpha}, B_{\alpha})$. 
So $\bar{a}_{\beta}$ realizes 
$p(\bar{x},\bar{c}_{\alpha},\bar{b})$.

Let $\bar{d}_{\alpha}=\bar{c}_{\alpha} \bar{a}_{\alpha}$
and let 
$q(\bar{x}_1, \bar{y}_1, \bar{x}_2, \bar{y}_2, \bar{b}):= 
p(\bar{x}_1, \bar{y}_2, \bar{b})$ 
(we may assume that $q$ is closed under finite 
conjunction). 
Then, above construction shows that
\begin{equation} \tag{*}
\bar{d}_{\alpha}\bar{d}_{\beta} \models  
q(\bar{x}_1, \bar{y}_1, \bar{x}_2, \bar{y}_2, \bar{b})
\quad
\text{ if and only if } 
\quad
\alpha < \beta < \mu,
\end{equation}
i.e. we we have an order of length $\mu$ witnessed
 by the type $q$.

We use (*) to obtain an order of 
length $\lambda$ witnessed by a formula as follows.
On the one hand,  (*) implies that for any 
$\phi(\bar{x}_1, \bar{x}_2, \bar{y}_1,\bar{y}_2, \bar{c}) 
\in q$, 
the following holds:
\begin{equation} \tag{**}
\models \phi[\bar{d}_\alpha, \bar{d}_\beta, \bar{b}]
\qquad 
\text{whenever $\alpha < \beta$}.
\end{equation}
On the other hand, if $\alpha \geq \beta$, by (*) again, there is 
$\phi_{\alpha, \beta}
(\bar{x}_1, \bar{x}_2, \bar{y}_1,\bar{y}_2, \bar{b}) 
\in q$, 
such that
$\models 
\neg \phi_{\alpha, \beta}[\bar{d}_\alpha, \bar{d}_\beta, \bar{b}]$.
Hence, by the Erd\"{o}s-Rado Theorem, since $|q| \leq |T|$, 
we can find $S \subseteq \mu$
of cardinality $\lambda$ and 
$\phi
(\bar{x}_1, \bar{x}_2, \bar{y}_1,\bar{y}_2, \bar{b}) 
\in q$, 
such
that 
\begin{equation} \tag{***}
\models \neg \phi[\bar{d}_\alpha, \bar{d}_\beta, \bar{b}]
\qquad
\text{whenever $\alpha \geq \beta$}, 
\quad
\alpha, \beta \in S.
\end{equation}
Therefore, (**) and (***) together show that we can find
an order of length $\lambda$, which is the desired contradiction.
\epf

We close this section by gathering together
the properties of the forking symbol.
They are stated with the names of the first-order 
forking properties to which they
correspond.

\begin{theorem} \label{forkingproperties} 
\begin{enumerate}
\item 
\textup{(Definition)} 
$A \nonfork_B C$ if and only if $A \nonfork_B B \cup C$.
\item 
\textup{(Existence)}
$A \nonfork_B B$
\item
\textup{($\kappa(D)=\aleph_0$)}
For all $\bar{a}$ and $C$, there is a finite $B \subseteq C$ 
such that $\bar{a} \nonfork_B C$.

\item 
\textup{(Invariance under automorphisms)} 
Let $f \in \aut{}$.
\[
A \nonfork_B C 
\qquad
\text{if and only if}
\qquad
 f(A) \nonfork_{f(B)} f(C).
\]
\item 
\textup{(Finite Character)} \label{finitechar}
\[
A \nonfork_B C 
\qquad
\text{if and only if}
\qquad
A' \nonfork_B C',
\] 
for every finite $A' \subseteq A$, and finite $C' \subseteq C$ .
\item 
\textup{(Monotonicity)} \label{mono} Suppose $A'$
and $C'$ contain $A$ and $C$ respectively and
that $B'$ is a subset of $B$. Then 
\[
A \nonfork_B C 
\qquad
\text{implies}
\qquad
 A' \nonfork_{B'} C'.
\]
\item 
\textup{(Transitivity)} If $B \subseteq C \subseteq D$, then
\[
A \nonfork_B C 
\quad
\text{and}
\quad
 A \nonfork_C D 
\qquad
\text{if and only if}
\qquad
A \nonfork_B D.
\]
\item 
\textup{(Symmetry)} \label{symm} Let $M$ is a $(D,\aleph_0)$-homogeneous model.
\[
A \nonfork_M C 
\qquad
\text{if and only if}
\qquad
 C \nonfork_M A.
\]
\item 
\textup{(Extension)}
Let $M$ be a $(D,\aleph_0)$-homogeneous model.
For every $A, C$ there exists $A'$ such that 
\[
\tp(A, M)=\tp(A',M) 
\quad
\text{ and }
\quad
 A' \nonfork_M C.
\]
\item
\textup{(Uniqueness)}
Let $M$ be a $(D, \aleph_0)$-homogeneous model.
If $A, A'$ satisfy 
\[
\tp(A,M)=\tp(A',M) 
\qquad
\text{and both}
\qquad
A \nonfork_M C 
\quad
\text{and}
\quad
A' \nonfork_M C
\]
then $\tp(A, MC)=\tp(A',MC)$.
\end{enumerate} 
\end{theorem}

\bpf 
\begin{enumerate}
\item 
This is just by Definition \ref{fork}.
\item
Immediate from Definition \ref{fork}.
\item 
By Finite Character of the rank and Definition \ref{fork}.
\item
Follows from Lemma \ref{rankproperties} \ref{rankinvariance}.
\item 
Immediate by finite definition and finite character of the rank.
\item 
Assume $C \fork_M A$. 
Then, by Finite Character, there is $\bar{c} \in C$, 
such that $R[\tp(\bar{c},M)]<R[\tp(\bar{c},M)]$. 
Also by Finite Character , there exists $\bar{a} \in A$ such that 
$R[\tp(\bar{c},M \cup \bar{a})]=R[\tp(\bar{c},M)]$. 
Hence $\bar{c} \fork_M \bar{a}$.
But, by Corollary \ref{stationary}, both $\tp(\bar{a},M)$ 
and $\tp(\bar{c},M)$ are stationary, 
so by Theorem \ref{symmetry} we must have $\bar{a} \fork_M \bar{c}$.
By Finite Character, this shows that $A \fork_M C$.
\item 
Let $\bar{a} \in A$. 
Then, by Finite Character, $\bar{a} \nonfork_B C$,
and $\bar{a} \nonfork_C D$, so by Definition \ref{fork}
$R[\tp(\bar{a},C)]=R[\tp(\bar{a}, B)]$ and 
$R[\tp(\bar{a},D)]=R[\tp(\bar{a},C)]$.
Thus $R[\tp(\bar{a}, B)] =R[\tp(\bar{a},D)]$, 
so $\bar{a} \nonfork_B D$.
Hence, by Finite Character, we must have $A \nonfork_B D$.
The converse is just by Monotonicity.
\item 
Immediate by Theorem \ref{symmetry} and Corollary \ref{stationary}.
\item 
Follows from Corollary \ref{stationary} and Definition \ref{fork}.
\item 
Follows from Corollary \ref{stationary} and Definition \ref{fork}.
\end{enumerate}
\epf

\section{Regular and Minimal types}

In this section, we prove the existence of various pregeometries for totally
transcendental diagrams. First, we make the following definition (a similar
definition appears in \cite{Sh:4}).

\begin{definition} \label{defbig}
\begin{enumerate}
\item 
Let $\bar{a}$ be in $M$ and
$q(\bar{x},\bar{a})$ be a type.
We say that $q(\bar{x},\bar{a})$ is \emph{big for $M$ }
if $q(\bar{x},\bar{a})$ is realized outside $M$;

\item 
We say that $q(\bar{x},\bar{a})$ is\emph{ big} if
$q(\bar{x},\bar{a})$ is big for any $M$ containing $\bar{a}$;

\item 
A type $q \in S_D(A)$ is \emph{big (for M)} if 
$q \restriction \bar{a}$ is big (for $M$) for every $\bar{a} \in A$.
\end{enumerate}
\end{definition}

In presence of the compactness theorem, big types are the same 
as non-algebraic types. 
Even in the general case, we have a nice characterization of bigness
when the types are stationary.

\begin{lemma} \label{bigone} Let $q \in S_D(A)$ be stationary. 
The following conditions are equivalent:
\begin{enumerate}
\item $q$ is big for some $(D,\aleph_0)$-homogeneous $M$ containing $A$;
\item $R[q] \geq 1$;
\item $q$ is big.
\end{enumerate} 
\end{lemma}
\bpf 
(1) $\Rightarrow$ (2): Since $M$ is $(D, \aleph_0)$-homogeneous, 
by Lemma \ref{homo},
$R[q]=R_M[q]$, so it is enough to show $R_M[q] \geq 1$.
Let $\bar{a} \in A$ be such that $R_M[q] = R_M [q \restriction \bar{a}]$.
Since $q \restriction \bar{a}$ is big for $M$,
there exists $\bar{c} \not \in M$
realizing $q \restriction \bar{a}$. 
Also, since $M$ is $(D,\aleph_0)$-homogeneous, 
there is $\bar{c}' \in M$ realizing 
$q \restriction \bar{a}$. 
Hence 
\[
R_M[(q \restriction \bar{a}) \cup \{\bar{x} = \bar{c}'\}] \geq 0
\quad
\text{ and }
\quad
R_M[(q \restriction \bar{a}) \cup \{ \bar{x} \not 
= \bar{c}' \}] \geq 0.
\] 
Moreover, for every $\bar{b} \in M$, 
$(q \restriction \bar{a}) \cup \tp(\bar{c},\bar{b})$ 
is realized by $\bar{c}$,
and so 
\[
R_M[(q \restriction \bar{a}) \cup \tp(\bar{c},\bar{b})]\geq 0,
\]
and $\tp(\bar{c},\bar{b}) \in S_D(\bar{b})$. 
This shows that $R_M[q \restriction \bar{a} ] \geq 1$.

(2) $\Rightarrow$ (3): Suppose $q$ is stationary,  $R[q] \geq 1$ and 
$M$ containing $\bar{a}$ are given. 
By taking a larger $M$ if necessary, we may assume that 
$M$ is $(D, \aleph_0)$-homogeneous.
Since $q$ is stationary, there exists $q_M \in S_D(M)$, such that
$R[q_M]=R[q]  \geq 1$.
Let $\bar{c}$ realize $q_M$. 
If $\bar{c} \in M$, then $\{ \, x = \bar{c}\,\} \in q_M$, 
so 
\[
0=R[\bar{x}=\bar{c}] \geq R[q_M] \geq 1,
\] 
which is  a contradiction.
Hence $\bar{c} \not \in M$, so $q$ is big for $M$.

(3) $\Rightarrow$ (1): Clear by definition.
\epf

\begin{definition} Let $p \in S_D(A)$ be a big, stationary type.
\begin{enumerate}
\item
We say that $p$ is\emph{ regular for $M$} if $A \subseteq M$
and for every $B \subseteq M$ we have
\[
\bar{a}  \nonfork_A B 
\text{ and }
\bar{b} \fork_A B
\quad
\text{imply}
\quad
\bar{a} \nonfork_A B \cup \bar{b},
\qquad 
\text{for all $\bar{a}, \bar{b} \in p(M)$.}
\]
\item We say that $p$ is \emph{regular} if $p$ is
regular for $\mathfrak{C}$.
\end{enumerate}
\end{definition}

\begin{lemma} Let $p \in S_D(A)$ be a big, 
stationary type based on $\bar{c} \in A$.
If $p \restriction \bar{c}$ is regular, then $p$ is regular.
\end{lemma}
\bpf First notice that stationarity and bigness are preserved 
(bigness is the content 
of Lemma \ref{bigone}).
Suppose $p$ is not regular. 
We will show that $p \restriction \bar{c}$ is not regular.
Let $\bar{a}, \bar{b} \models p$ and $B$ be such that
\[
\bar{a}  \nonfork_A B,
\quad
\bar{b} \fork_A B
\quad
\text{and yet}
\quad
\bar{a} \fork_A B \cup \bar{b}.
\]
Therefore $\tp(\bar{a}, A \cup B) = p_{A \cup B}$
and so by choice of $\bar{c}$
we have  $\tp(\bar{a}, A \cup B) = (p \restriction{c})_{A \cup B}$,
i.e. $\bar{a} \nonfork_{\bar{c}} A \cup B$.
Now since $R[p] =R[ p\restriction \bar{c}]$,
\[
R[\tp(\bar{b}, A \cup B)] < R[\tp(\bar{b}, A)]
\quad
\text{implies}
\quad
R[\tp(\bar{b}, A \cup B)] < R[p\restriction \bar{c}],
\]
i.e. $\bar{b} \fork_{\bar{c}} A \cup B$.
We show similarly that $\bar{a} \fork_{\bar{c}} A \cup B \cup \bar{b}$,
which shows that $p \restriction \bar{c}$ is not regular.
\epf

\begin{remark} If $p(\bar{x}, \bar{a})$ is regular 
and $\bar{a}' \in M$ is such that
$\tp(\bar{a}, \empty)=\tp(\bar{a}', \empty)$,
then $p(\bar{x}, \bar{a}')$ is regular.
\end{remark}

\begin{definition}  \label{pregeometry} 
Let $p \in S_D(B)$, $B \subseteq M$ and $W = p(M)-B \not = \empty$.
Define
\[
a \in cl(C)
\quad
\text{ if }
\quad
a \fork_B C,
\qquad 
\text{for $a \in W$ and $C \subseteq W$}.
\]
\end{definition}

\begin{theorem} Let 
$M$ be $(D, \aleph_0)$-homogeneous containing   
$B$ and $p \in S_D(B)$ be 
realized in $M$.
If $p$ is regular then $(W,cl)$ is a pregeometry.
\end{theorem}

\bpf We need to show that the four axioms of pregeometry hold 
(notice that $W \not = \empty$).
\begin{enumerate}
\item 
We show that for every $C \subseteq W$, $C \subseteq cl(C)$.

Let $c \in C$, then $\{x=c\} \in \tp(c, A \cup C)$, 
hence 
\[
R[\tp(c, B \cup C)]=0 < R[p],
\]
so $c \fork_B C$ and thus $c \in cl(C)$.
\item 
We show that if $c \in cl(C)$, there is $C' \subseteq C$ finite, 
such that $c \in cl(C')$.

Let $c \in cl(C)$. By Definition \ref{pregeometry} $c \fork_B C$ so 
by Theorem \ref{forkingproperties} \ref{finitechar}
there exists $C' \subseteq C$ finite, such that $c \fork_B C'$, 
hence $c \in cl(C')$.
\item 
We show that if $a \in cl(C)$ and $C \subseteq cl(E)$, then $a \in cl(E)$.

Write $C =\{c_i \mid i < \alpha \}$. Then $a \fork_B \{c_i \mid i < \alpha \}$.
Suppose $a \nonfork_B E$.
We show by induction on $i < \alpha$ that 
$a \nonfork_B E \cup \{ c_j \mid j<i\}$.
\begin{itemize}
\item
For $i=0$ this is the assumption and for $i$ a limit ordinal, this is true by 
Theorem \ref{forkingproperties} \ref{finitechar}.
\item
For the successor case, suppose it is true for
$i$. 
Then  $a \nonfork_B E \cup \{c_l \mid l< i\}$.
Since $C \subseteq cl(E)$, we have $c_i  \fork_B E$, 
so by Theorem \ref{forkingproperties} \ref{mono}
$c_i \fork_B E \cup \{c_l \mid l<i\}$.
Hence, since $p$ is regular, we must have 
$a \nonfork_B E \cup \{ c_l \mid l <i  \}\cup c_i$.
\end{itemize}
Thus $a \nonfork_B E\cup C$, and since $C \subseteq C \cup E$, we must
have $a \nonfork_B C$. 
Hence $a \not \in cl(C)$, which contradicts our assumption.
\item 
We show that if $c \in cl(Ca)-cl(C)$, then $a \in cl(Cc)$.

Since symmetry has been shown only for stationary types, 
this statement is not immediate from 
Theorem \ref{symmetry}.

Suppose that $c \fork_B Ca$ and $c \nonfork_B C$. 
Then $c \fork_C a$, since
\[
R[ \tp(c, B\cup Ca)] < R[\tp(c, B)] = R[\tp(c, B \cup C)].
\] 
Therefore $c$ realizes $p_{B \cup C}$, so $\tp(c, B \cup C)$ is stationary.
If $a \fork_B C$, then by Theorem \ref{forkingproperties} \ref{mono} 
we must have $a \fork_B Cc$, 
and we are done.

Otherwise, $a \nonfork_B C$. 
Hence $a$ realizes $p_{B \cup C}$ and so $\tp(a, B \cup C)$ is
stationary. 
Therefore by Theorem \ref{symmetry} we must have 
$a \fork_C c$, a contradiction.
Hence by Theorem \ref{forkingproperties} \ref{mono}, 
we have $a \fork_B Cc$, i.e. $a \in cl(Cc)$.
\end{enumerate} 
\epf

We now show the connection between independent sets in the pregeometries,
averages and stationarity.

\begin{lemma} \label{baseaverage}Let $p(\bar{x}, \bar{c})$ be regular.
Suppose $I$ is infinite and independent in $p(\mathfrak{C}, \bar{c})$.
Then $I$ is indiscernible and for every $B$ containing $\bar{c}$ we have
$p_B=\av(I,B)$.
\end{lemma}
\bpf Write $I=\{\bar{a}_i \mid i < \alpha \}$.
Then since $I$ is independent, $\bar{a}_{i+1} \models p_{A_i}$, 
where $A_i= \bar{c} \cup \{\bar{a}_j \mid j < i \}$.
Thus $I$ is a Morley sequence based on $p$, so the result 
follows from Lemmas \ref{morley}
and \ref{averagemorley}.
\epf

Now we turn to existence. In order to do this, we need a lemma.

\begin{lemma} \label{regularinM} Let $M$ 
be $(D,\aleph_0)$-homogeneous,
and $p(\bar{x}, \bar{c})$ over $M$
be big and stationary.
Then $p(\bar{x}, \bar{c})$ is regular if and only if
$p(\bar{x},\bar{c})$ is regular for $M$.
\end{lemma}
\bpf If $p(\bar{x},\bar{c})$ is regular, then
$p(\bar{x},\bar{c})$ is clearly regular for $M$.
Suppose $p(\bar{x}, \bar{c})$ is not regular.
Then there are $B \subseteq \mathfrak{C}$, and  $\bar{a}$, $\bar{b}$ 
realizing $p(\bar{x}, \bar{c})$, such that
\[
\bar{a} \nonfork_{\bar{c}} B, 
\quad
\bar{b} \fork_{\bar{c}} B,
\quad
\text{and}
\quad
\bar{a} \fork_{\bar{c}} B \bar{b}.
\]
First, we may assume that $B$ is finite: 
choose $B' \subseteq B$ such that
\[
R[\tp(\bar{a}, B' \cup \bar{c} \bar{b})]=
R[\tp(\bar{a}, B \cup \bar{c} \bar{b})]
\] 
and then choose $B'' \subseteq B$ finite, 
such that $\bar{b} \not \models p_B \restriction B''$.
Hence, for $B_0=B' \cup B'' \subseteq B$, 
we have  
\[
\bar{a} \nonfork_{\bar{c}} B_0, 
\quad
\bar{b} \fork_{\bar{c}} B_0,
\quad
\text{and}
\quad
\bar{a} \fork_{\bar{c}} B_0 \bar{b}.
\]
Now, since $M$ is $(D, \aleph_0)$-homogeneous and $\bar{c} \in M$, 
we can find $B_1, \bar{a}_1$ and $\bar{b}_1$ inside $M$ such that
$\tp (B_0 \bar{a} \bar{b} , \bar{c})=\tp(B_1 \bar{a}_1 \bar{b}_1 , \bar{c})$. 
Therefore, by invariance we have:
\[
\bar{a} \nonfork_{\bar{c}} B_1, 
\quad
\bar{b} \fork_{\bar{c}} B_1,
\quad
\text{and}
\quad
\bar{a} \fork_{\bar{c}} B_1 \bar{b}.
\]
This shows that $p$ is not regular for $M$.
\epf

The following argument for the existence of regular types is similar to 
Claim V.3.5. of \cite{Sh:a}.
However, since our basic definitions are 
different, we provide a proof.

\begin{theorem} [Existence of regular types]\label{regular}
Let $M \subseteq N$ be 
$(D, \aleph_0)$-homogeneous.
If $M \not = N$, then there exists $p(x, \bar{a})$ regular, realized in $N -M$.
In fact, if $p(x, \bar{a})$ is big and stationary, 
and has minimal rank among all
big, stationary types over $M$ realized in $N-M$, 
then $p(x,\bar{a})$ is regular.
\end{theorem} 
\bpf The first statement follows from the second.
To prove the second statement, we first choose 
$c' \in N-M$, be such that $\tp(c',M)$ has minimal rank
among all types
over $M$ realized in $N-M$, say $R[\tp(c',M)]=\alpha$.
We then choose $\bar{a} \in M$ such that 
$R[\tp(c',M)]=R[\tp(c', \bar{a})]=\alpha$.
Write $\tp(c', \bar{a})=p(x, \bar{a})$ and notice that $p$ is stationary
and big for $M$, hence big, by Lemma \ref{bigone}.

By the previous lemma, to show that $p(x,\bar{a})$
is regular, it is equivalent to show that 
$p(x,\bar{a})$ is regular for $M$. For this, 
let $a, b \in p(M)$ and $B \subseteq M$ 
such that 
\[
a \nonfork_{\bar{a}} B
\quad 
\text{ and }
\quad
b \fork_{\bar{a}} B.
\]
We must show that $a \nonfork_{\bar{a}} Bb$. 
Suppose, by way of contradiction that this is not the case.
Then, by definition, we have  $R[\tp(a,B\bar{a}b)]<\alpha$. 
We now choose $\bar{c}, \bar{d} \in B$ such that 
\[
R[\tp(a,B\bar{a}b)]=R[\tp(a,\bar{c}\bar{a}b)]<\alpha 
\quad
\text{ and }
\quad
R[\tp(b,B\bar{a})]=R[\tp(b,\bar{d}\bar{a})]<\alpha.
\]
Since $N$ is $(D,\aleph_0)$-homogeneous and 
$c',a, b, \bar{a}, \bar{c}, \bar{d} \in N$,
there is $b' \in N$ such that 
$\tp(ab, \bar{a}\bar{c}\bar{d})=
\tp(a'b', \bar{a}\bar{c}\bar{d})$.
Now, $\tp(b', \bar{a}\bar{d})=\tp(b', \bar{a}\bar{d})$, so 
\[
R[\tp(b',M)] \leq R[\tp(b', \bar{a}\bar{d})]=
R[\tp(b, \bar{a}\bar{d})] < \alpha.
\]
By minimality of $\alpha$, we must have $b' \in M$.
This implies that $R[\tp(a', M)] \leq R[\tp(a',\bar{c}\bar{a}b')]$, so
$R[\tp(a',\bar{c}\bar{a}b')]=\alpha$.
Now there is $f \in \aut{}$ such that $f(a')=a$, $f(b')=b$ and 
$f \restriction \bar{c}\bar{a} =
id_{\bar{c}\bar{a}}$, by choice of $b'$. 
Hence, by property of the rank
\[
\alpha = R[\tp(a',\bar{c}\bar{a}b')]= 
R[ f(\tp(a',\bar{c}\bar{a}b'))]=  
R[\tp(a,\bar{c}\bar{a}b)]<\alpha,
\] 
which is a contradiction.
Hence $a \nonfork_{\bar{a}} Bb$, so that $p(x, \bar{a})$
is regular.
\epf

By observing what happens when $N = \mathfrak{C}$ in above theorem,
 one discovers more
concrete regular types. For this, we make the following definition.
A similar definition in the context of $L_{\omega_1\omega}(Q)$ appears 
in the last section of \cite{Sh:4}. An illustration of why this definition
is natural can be found in the proof of Lemma \ref{saturation}.
In presence of the compactness theorem, S-minimal 
is the same as strongly minimal.

\begin{definition}
\begin{enumerate}
\item 
A big, stationary type $q(\bar{x}, \bar{a})$ over $M$
is said to be \emph{S-minimal for $M$} if for any 
$\theta(\bar{x},\bar{b})$ over $M$ 
not both 
$q(\bar{x},\bar{a}) \cup \theta(\bar{x},\bar{b})$ and
 $q(\bar{x},\bar{a}) \cup \neg \theta(\bar{x},\bar{b})$ 
are big for $M$.
\item 
A big, stationary type $q(\bar{x}, \bar{a})$
is said to be \emph{S-minimal} if $q(\bar{x}, \bar{a})$ is S-minimal for
for every $M$ containing $\bar{a}$.
\item  
If $q \in S_D(A)$ is big and stationary, we say that $q$ is \emph{S-minimal}
if $q \restriction \bar{a}$ is S-minimal for some $\bar{a}$.
\end{enumerate}
\end{definition}

\begin{remark}\label{Sminimal}
\begin{enumerate}
\item 
Let $q(\bar{x}, \bar{c})$ be S-minimal for the
$(D,\aleph_0)$-homogeneous model $M$.
Let $W= q(M, \bar{c})$
and for $a \in W$ and $B \subseteq W$ define
\[
a \in cl(B) 
\quad
\text{ if }
\quad
\tp(a, B \cup \bar{c}) \text{ is not big (for $M$)}.
\]
Then it can be shown directly from the assumption that $D$
is totally transcendental, that $(W, cl)$ is a pregeometry.
\item \label{Sminimal1} 
If $M$ is $(D, \aleph_0)$-homogeneous and $q(x, \bar{c})$
has minimal rank among all big, stationary $q(x, \bar{c})$
over $M$, then 
the previous theorem shows that $q$ is regular.
But $q$ is also S-minimal for $M$.
As a matter of fact, if $a \nonfork_{\bar{c}} B$, then 
$R[\tp(a, B \cup \bar{c})]= R[q(\bar{x}, \bar{c})] \geq 1$ and
$\tp(a, B \cup \bar{c})$ is stationary, so $\tp(a, b \cup \bar{c})$
is big, so $a \not \in cl(B)$. 
Conversely, if $a \fork_{\bar{c}} B$, then 
$R[\tp(a, B\bar{c})] < R[q(x, \bar{c}]$. 
But if $\tp(a, B \cup \bar{c})$ was big, then we could find
$a' \not \in M$ such that 
$\tp(a', B \cup \bar{c})=\tp(a, B \cup \bar{c})$,
so 
\[
R[\tp(a', M)] \leq R[\tp(a', B \cup \bar{c})] = 
R[\tp(a, B \cup \bar{c})] < R[q(x, \bar{c})],
\]
contradicting the minimality of $R[q(x, \bar{c})]$.
Hence $\tp(a, B \cup \bar{c})$ is not big, and so $a \in cl(B)$.
In other words, both pregeometries coincide.
\item \label{Sminimal2}
Using the results that we have proven so far, it is not difficult
to show that if $M, N$ are $(D, \aleph_0)$-homogeneous, 
and $q(x, \bar{c})$ has minimal rank among all big, stationary types
over $M$ and $\bar{c}' \in N$ such that
$\tp(\bar{c}, \empty)=\tp(\bar{c}',\empty)$, then
$q(x, \bar{c}')$ has minimal rank among all big, stationary
types over $N$, hence if $q(x,\bar{c}')$ is S-minimal for $N$.
\end{enumerate}
\end{remark}

In the light of these remarks, we will make the following definition.
\begin{definition} Let $M$ be $(D,\aleph_0)$-homogeneous.
A big, stationary type $q(\bar{x}, \bar{c})$ with $\bar{c} \in M$
is called\emph{ minimal} if
$q(\bar{x}, \bar{c})$ has minimal rank among all  big, stationary
types over $M$.
\end{definition}

We close this section by summarizing above remark in the following theorem.

\begin{theorem} \label{minimalregular}
\begin{enumerate}
\item 
For any $(D, \aleph_0)$-homogeneous
model, there exists a minimal $q(x, \bar{c})$ with $\bar{c} \in M$.
\item 
Minimal types are regular and moreover for every $A$ containing  
$\bar{c}$, every set $B$
and $a \models q_A$ we have
\[
\tp(a, A \cup B) 
\text{ is big }
\quad
\text{ if and only if }
\quad
a \nonfork_A B.
\]
\end{enumerate}
\end{theorem}
\bpf The first item is clear by definition.
The second follows by Theorem \ref{regular}, and 
Remark \ref{Sminimal} \ref{Sminimal1}
and \ref{Sminimal2}.
\epf

\section{Applications}

In this section, we give a few applications of our concepts.
The rank is especially useful to study the class of $(D,\aleph_0)$-homogeneous
models  of a totally transcendental $D$. 
In the first subsection, we start with the existence of prime models.

\subsection{Prime models}

We give definitions from \cite{Sh:1} in
more modern terminology.

\begin{definition} 
\begin{enumerate}
\item
We say that $p\in S_D(A)$ is \emph{$D^s_{\lambda}$-isolated over 
$B \subseteq A$}, $|B|<\lambda$, 
if for any $q \in S_D(A)$ extending $p \restriction B$, we have $q=p$.

\item
We say that $p\in S_D(A)$ is\emph{ $D^s_{\lambda}$-isolated} if there is 
$B \subseteq A$, $|B|<\lambda$, 
such that $p$ is $D^s_{\lambda}$-isolated over $B$.
\end{enumerate}  
\end{definition}
The following are verifications of Axioms X.1 and XI.1 
from Chapter IV of \cite{Sh:a}. 

\begin{theorem} [X.1] Let 
$A \subseteq \mathfrak{C}$ and 
$\mu \geq \aleph_0$.
Every $\phi(\bar{x},\bar{a})$ over $A$ realized in $\mathfrak{C}$
can be extended to a $D^s_{\mu}$-isolated type $p \in S_D(A)$.
\end{theorem}
\bpf It is enough to show the result for $\mu = \aleph_0$.

Since $\mathfrak{C} \models \exists \bar{x} \phi[\bar{x},\bar{a}]$, 
there exists
$\bar{c} \in \mathfrak{C}$ such that 
$\mathfrak{C} \models \phi[\bar{c},\bar{a}]$.
Thus there exists is $p\in S_D(A)$, namely $\tp(\bar{c},A)$, containing 
$\phi(\bar{x},\bar{a})$.
Since $D$ is totally transcendental and $A \subseteq \mathfrak{C}$ 
we must have $R_A[p]<\infinity$.
Among all those $p\in S_D(A)$ containing $\phi(\bar{x},\bar{a})$ 
choose one with minimal rank.
Say $R_A[p]=\alpha\geq 0$.

We claim that $p$ is $D^s_{\aleph_0}$-isolated.
First, there is $\bar{b} \in A$ such that 
$R_A[p]=R_A[p \restriction \bar{b}]$. 
We may assume that 
$p \restriction \bar{b}$ contains $\phi(\bar{x},\bar{a})$ 
by Lemma \ref{rankproperties} \ref{subset}.
Suppose that there is $q \in S_D(A)$, $q \not = p$, 
such that $q$ extends $p \restriction \bar{b}$.
Then $R_A[q]\geq \alpha$ by choice of $p$ 
(since $q$ contains $\phi(\bar{x},\bar{a})$).
Now, choose $\psi(\bar{x},\bar{c})$ with $\bar{c} \in A$ such that 
$\psi(\bar{x},\bar{c}) \in p$ and 
$\neg \psi(\bar{x},\bar{c}) \in q$.
Then 
since 
$(p \restriction \bar{b}) \cup \psi(\bar{x},\bar{c}) \subseteq p$,
by Lemma \ref{rankproperties} \ref{subset} we have
\[
R_A[(p \restriction \bar{b}) \cup \psi(\bar{x},\bar{c})] \geq R_A[p]
\geq \alpha.
\] 
Similarly 
\[
R_A[(p \restriction \bar{b}) \cup \neg \psi(\bar{x},\bar{c})] 
\geq R_A[q]\geq \alpha.
\]
Now, given any $\bar{d} \in A$, 
$R_A[p \restriction \bar{b}\cup\bar{d}]\geq \alpha$
(again by Lemma \ref{rankproperties} \ref{subset}). 
Since $p \in S_D(A)$, necessarily if we write 
$p \restriction \bar{d}= p(\bar{x},\bar{d})$, then we have 
$p(\bar{x},\bar{y}) \in D$ 
(since $p(\bar{x},\bar{d}) \in S_D(\bar{d})$).
Hence since $p \restriction \bar{b}\cup\bar{d} 
\vdash p \restriction {b} \cup p(\bar{x},\bar{d})$)
we have
\[
R_A[(p \restriction {b}) \cup p(\bar{x},\bar{d})]\geq 
R_A[p \restriction \bar{b}\cup\bar{d}]\geq \alpha.
\] 
But this shows that $R_A[p \restriction \bar{b}]\geq \alpha +1$, 
a contradiction.

Hence $p$ is the only extension of $p \restriction {b}$, 
so $p$ is $D^s_{\aleph_0}$-isolated.  
\epf

\begin{theorem} [XI.1] Let $\mu$ be infinite and $B \subseteq A$.
Every $D^s_{\mu}$-isolated $r~\in~S_D(B)$ can be extended
to a $D^s_{\mu}$-isolated type $p \in S_D(A)$.
\end{theorem}

\bpf Since $\mathfrak{C}$ is $(D,\chi)$-homogeneous, there exists 
$\bar{c} \in \mathfrak{C}$ realizing $r$. 
Hence there is $p\in S_D(A)$ extending $r$, namely $\tp(\bar{c},A)$.
Since $D$ is totally transcendental and $A \subseteq \mathfrak{C}$ 
we must have $R_A[p]<\infinity$.
Among all those $p\in S_D(A)$ extending $r$ choose one with minimal rank.
Say $R_A[p]=\alpha\geq 0$.

We claim that $p$ is $D^s_{\mu}$-isolated.
First, there is $\bar{b} \in A$ such that 
$R_A[p]=R_A[p \restriction \bar{b}]$. 
Also, since $r$ is $D^s_\mu$-isolated, there is $C \subseteq B$,
$|C| < \mu$ such that $r \restriction C$ isolates $r$.
We may assume that $R_A[r] = R_A[r \restriction C]$, 
by Lemma \ref{rankproperties} \ref{finite}.
We claim that $(r \restriction C) \cup (p \restriction \bar{b})$ isolates $p$.
By contradiction, suppose that there is $q \in S_D(A)$ 
extending $(r \restriction C) \cup (p \restriction \bar{b})$
such that $q \not = p$.
Notice that $r \subseteq q$, since $r$ was isolated by $r \restriction C$,
 and hence $R_A[q]\geq R_A[p] =\alpha$ by choice of $p$.
Now, choose $\psi(\bar{x},\bar{a})$ with $\bar{a} \in A$ such that 
$\psi(\bar{x},\bar{a}) \in p$ and 
$\neg \psi(\bar{x},\bar{a}) \in q$.
By Lemma \ref{rankproperties} \ref{subset} 
(since 
$(p \restriction \bar{b}) \cup \psi(\bar{x},\bar{c}) \subseteq p$), 
we must have 
\[
R_A[(p \restriction \bar{b}) \cup \psi(\bar{x},\bar{c})] \geq R_A[p]
=\alpha.
\]
Similarly 
\[
R_A[(p \restriction \bar{b}) \cup \neg \psi(\bar{x},\bar{c})] 
\geq R_A[q]\geq \alpha.
\]
Now, given any $\bar{d} \in A$ we have that 
$R_A[p \restriction \bar{b}\cup\bar{d}]\geq \alpha$
(again by Lemma \ref{rankproperties} \ref{subset}). 
Since $p \in S_D(A)$, necessarily if we write 
$p \restriction \bar{d}= p(\bar{x},\bar{d})$, then we have 
$p(\bar{x},\bar{y}) \in D$ 
(since $p(\bar{x},\bar{d}) \in S_D(\bar{d})$).
Hence 
\[
R_A[( p \restriction {b}) \cup p(\bar{x},\bar{d})]\geq 
R_A[p \restriction \bar{b}\cup\bar{d}]\geq \alpha,
\]
since $p \restriction \bar{b}\cup\bar{d} \vdash 
(p \restriction {b}) \cup p(\bar{x},\bar{d})$.
But this shows that $R_A[p \restriction \bar{b}]\geq \alpha +1$, 
a contradiction.

Hence $p$ is the only extension of 
$(r \restriction C) \cup (p \restriction {b})$, 
so $p$ is $D^s_{\mu}$-isolated.  
\epf

Following Chapter IV of \cite{Sh:a}, we set:
\begin{definition}
\begin{enumerate}
\item 
We say that $\mathcal{C}$=$\{\langle a_i,A_i,B_i \rangle \mid i<\alpha\}$ is a 
\emph{$(D, \lambda)$-construction of $C$ over $A$} if 
\begin{enumerate}
\item
$C= A \cup \Union \{ a_i \mid i<\alpha\}$;
\item
$B_i \subseteq A_i$, $|B_i| < \lambda$, where
$A_i= A \cup \Union \{ a_j \mid j<i \}$; 
\item 
$\tp(a_i,A_i) \in S_D(A_i)$ is
$D^s_\lambda$-isolated over $B_i$.
\end{enumerate}

\item 
We say that $M$  is\emph{ $D_{\lambda}^s$-constructible over $A$} if 
there is a 
$(D,\lambda)$-construction for $M$ over $A$.
\item
 We say that $M$ is \emph{$D_{\lambda}^s$-primary over $A$}, if $M$ is 
$D_{\lambda}^s$-constructible over $A$ and $M$ is $(D,\lambda)$-homogeneous.
\item 
We say that $M$ is \emph{$D_{\lambda}^s$-prime over $A$} if 
\begin{enumerate}
\item 
$M$ is $(D,\lambda)$-homogeneous and 
\item 
if $N$ is $(D,\lambda)$-homogeneous and $A \subseteq N$, 
then there is $f: N \rightarrow M$ elementary 
such that $f \restriction A =id_A$.
\end{enumerate}
\item 
We say that $M$ is \emph{$D^s_{\lambda}$-minimal over $A$}, if $M$
is $D^s_{\lambda}$-prime over $A$ and for every $(D,\lambda)$-homogeneous model
$N$, 
if $A \subseteq N \subseteq M$, then $M=N$.
\end{enumerate}
\end{definition}

\begin{remark} We use the same notation as in \cite{Sh:a}, 
except that we replace ${\bf F}$ by
$D$ to make it explicit that we deal exclusively with $D$-types 
(or equivalently, types realized in
$\mathfrak{C}$). 
In particular, for example if $M$ is $D_{\aleph_0}^s$-primary over $A$, then
$M$ is $D_{\aleph_0}^s$-prime over $A$.
\end{remark}

\begin{theorem} [Existence of prime models] \label{prime}
Let $D$ be totally transcendental.
Then for all $A \subseteq \mathfrak{C}$ and infinite $\mu$  
there is a
$D_{\mu}^s$-primary model $M$ over $A$ of cardinality $|A|+|T|+|D|+\mu$.
Moreover, $M$ is $D^s_\mu$-prime over $A$.
\end{theorem}
\bpf See page 175 of \cite{Sh:a} and notice that we just established $X.1$ 
and $XI.1$. 
Observe that in the construction, each new element realizes a $D$-type, 
so that the resulting model
is indeed a $D$-model. The optimal bound on the cardinality follows from
Theorem \ref{stable}. The second sentence follows automatically.
\epf

\begin{remark} A similar theorem, with a stronger assumption 
($D$ is $\aleph_0$-stable) and without
the bound on the cardinality appears in \cite{Sh:1}.
Note that $D^s_{\mu}$-primary, is called $(D,\mu, 1)$-prime there.
\end{remark}

Notice that this allows us to show how any type can be
decomposed into stationary and isolated types. 
A similar result appears in \cite{Sh:871}.

\begin{lemma} \label{nice} 
Let $p \in S_D(A)$ and suppose $\bar{a}$ realizes $p$.
Then  there is $\bar{b} \in \mathfrak{C}$ such that
\begin{enumerate}
\item $\tp(\bar{b},A)$ is $D^s_{\aleph_0}$-isolated;
\item $\tp(\bar{a}, A \bar{b})$ is stationary;
\item $R[\tp(\bar{a}, A \bar{b})]=R[\tp(\bar{a},\bar{b})]$.
\end{enumerate} 
Furthermore, $p$ does not split over a finite set.
\end{lemma}
\bpf Let $\bar{a} \models p$.
Let $M$ be $D^s_{\aleph_0}$-primary model over $A$. 
Then $\tp(\bar{a},M)$ is stationary since $M$ is $(D, \aleph_0)$-homogeneous,
 and there is $\bar{b} \in M$ finite, such
that $R[\tp(\bar{a},M)]=R[\tp(\bar{a}, \bar{b})]$.
Hence $R[\tp(\bar{a}, A \bar{b})]=R[\tp(\bar{a},\bar{b})]$ 
by Lemma \ref{rankproperties} \ref{subset}, and so
$\tp(\bar{a}, A \bar{b})$ is stationary. 
Also, $\tp(\bar{b},A)$ is $D^s_{\aleph_0}$-isolated,
since $M$ is $D^s_{\aleph_0}$-primary over $A$.

Finally, to see that $p$ does not split over a finite set, 
assume $\bar{a} \models p$,
$\tp(\bar{b},A)$ is $D^s_{\aleph_0}$-isolated,
$\tp(\bar{a}, A \bar{b})$ is stationary, and 
$R[\tp(\bar{a}, A \bar{b})]=R[\tp(\bar{a},\bar{b})]$. 
Then there is $C \subseteq A$ finite, 
such that $\tp(\bar{b},A)$ is $D^s_{\aleph_0}$-isolated over $C$.
Also, since $\tp(\bar{a}, A \bar{b})$ is stationary, 
it does not split over $\bar{b}$.
Now it is easy to see that $p$ does not split over $C$: 
otherwise there are $\bar{c}_l \in A$,
and $\phi(\bar{x}, \bar{y})$
such that $\tp(\bar{c}_1,C)= \tp(\bar{c}_2,C)$, 
$\bar{c}_l \in A$ for $l=1,2$, and 
$\models \phi[\bar{a}, \bar{c}_1]$ and
$\models \neg \phi[\bar{a}, \bar{c}_2]$.
But $\tp(\bar{b},A)$ does not split over $C$, and so 
$\tp(\bar{c}_1,\bar{b})= \tp(\bar{c}_2,\bar{b})$.
However, this contradicts the fact that $\tp(\bar{a}, A \bar{b})$ 
does not split over $\bar{b}$. All the conditions are satisfied.
\epf

This gives us an alternative and short proof that averages are well-defined,
and in fact, allows us to give short proofs of all the facts in 
Lemma \ref{kappa}.

\begin{lemma} Let $I$ be infinite and $A \subseteq \mathfrak{C}$. 
Then $\av(I,A) \in S_D(A)$ 
\end{lemma}
\bpf Completeness is clear. 
To see that $\av(I,A)$ is consistent, suppose that both 
$\phi(x, \bar{a})$ and $\neg \phi(x, \bar{a})$
are realized by infinitely many elements of $I$.
But $\tp(\bar{a}, I)$ does not split over a finite set 
$B \subseteq I$ by the previous lemma.
Hence, by choice of $\phi(x, \bar{a})$, we can find
$b,c \in I-B$ such that $\models \phi[b, \bar{a}]$ and 
$\models \neg \phi[c, \bar{a}]$.
This however, shows that 
$\tp(\bar{a}, I)$ splits over $B$, since $\tp(b,B)=\tp(c,B)$ by
indiscernibility of $I$ and both $\phi(b, \bar{y})$, 
$\neg \phi(c, \bar{y}) \in \tp(\bar{a}, I)$.
Now $\av(I,A) \in S_D(A)$ since we can extend $I$ to a $D$-set of 
indiscernible $J$
of cardinality $|A|^+$, and then some element of $J$ realizes $\av(I,A)$.
\epf

The following is a particular case of  Theorem \ref{union}. 
We include it here not just for completeness,
 but because the proof is different from
the proof of \ref{union} and very similar 
in the conceptual framework to the first-order case.

\begin{theorem} \label{transunion}
Let $D$ be totally transcendental.
If $\langle M_i \mid i < \alpha \rangle$ is an increasing chain of 
$(D,\mu)$-homogeneous models,
then $\Union_{i < \alpha} M_i$ is $(D, \mu)$-homogeneous ($\mu$ infinite).
\end{theorem}
\bpf Let $M = \Union_{i<\alpha} M_i$ and notice that
$M$ is $(D,\aleph_0)$-homogeneous.
Let $p \in S_D(A)$, $A \subseteq M$, $|A| < \mu$ and
 choose $q \in S_D(M)$ extending $p$.
Then, by Corollary \ref{stationary}, $q$ is stationary 
and there is $B \subseteq M$,
finite such that $q$ is based on $B$. 
Let $i < \alpha$, be such that $B \subseteq M_i$.
Since $M_i$ is $(D, \mu)$-homogeneous, 
there is $I =\{a_j \mid j< \mu \} \subseteq M_i$
a Morley sequence for $q_B$. 
Then, by Lemma \ref{averagemorley}, 
$q_{AB} =\av(I,A\cup B)$.
But $|I| > |A \cup B|$, so by Lemma \ref{kappa} 
there is $a_j \in I$ realizing $\av(I,A\cup B)$.
But $q_{AB} \supseteq p$, so
$p$ is realized in $M$.
This shows that $M$ is $(D,\mu)$-homogeneous.
\epf

\subsection{Categoricity}

We now focus on the structure of $(D,\aleph_0)$-homogeneous models.
Notice that when $D$ is the set of isolated types over the empty set or when
$D$ comes from a Scott sentence of $L_{\omega_1\omega}$, this class
coincides with the class of 
$D$-models. When $D=D(T)$, then $\mathcal{K}$ is the class of
$\aleph_0$-saturated models (of a totally transcendental
theory, in our case).

\begin{definition} Define 
\[
\mathcal{K}=\{\, M \mid M \text{ is }(D,\aleph_0)-\text{homogeneous }\,\}.
\]
\end{definition}

\begin{remark} We will say that $M \in \mathcal{K}$ is
 \emph{prime over $A$} 
or \emph{minimal over $A$},
when $M$ is $D^s_{\aleph_0}$-prime over $A$ 
or $D^s_{\aleph_0}$-minimal over $A$ respectively.
\end{remark}

By analogy with the first-order case, we set the following definition.

\begin{definition}  Let $D$ be totally transcendental.
We say that $D$ is\emph{ unidimensional} if 
for every pair of models $M \subseteq N$ in $\mathcal{K}$ 
and minimal type $q(x,\bar{a})$ minimal over $M$, 
\[
q(M, \bar{a})=q(N, \bar{a})
\quad
 \text{ implies }
\quad
 M=N.
\]
\end{definition}

Unidimensionality for a totally transcendental diagram
 $D$ turns out to be a weak
dividing line. When it fails, we can construct non-isomorphic models,
like in the next theorem (this justifies the name), and when it holds 
we get a strong structural theorem
(see Theorem \ref{minimal}, which implies categoricity).
In fact, the conclusion of our next theorem is similar to (but stronger than) 
the conclusion of
Theorem 6.9 of \cite{Sh:1} (we prove it for every $\mu$, 
not just regular $\mu$,
and can obtain these models of cardinality exactly $\lambda$, 
not arbitrarily large).
The assumptions of Theorem 6.9 of \cite{Sh:1} are weaker 
and the proof considerably
longer. Actually, Corollary \ref{corollary3} 
makes the connection with Theorem 6.9
of \cite{Sh:1} clearer.

We first prove two technical lemmas which are 
similar to Lemma 3.4 and fact 3.2.1 from \cite{GrHa} respectively.
The proofs are straightforward generalizations and are presented here
for the sake of completeness.

\begin{lemma} \label{tech} Let 
$p, q \in S_D(M)$ and $M \subseteq N$ be in $\mathcal{K}$.
If
$a \nonfork_M b$
for every $a \models q$ and $b \models p$,
then 
$a \nonfork_N b$ for every $a \models q_N$ and $b \models p_N$.
\end{lemma}
\bpf Suppose not. Then there are $a \models p_N$ and 
$b \models q_N$ such that $a \fork_N b$.
Choose $E \subseteq N$ finite such that $a \fork_{ME} b$ and
$\tp(ab,N)$ is based on $E$. 
This is possible by Theorem \ref{forkingproperties} \ref{finitechar} and by
the fact that $\tp(ab,N)$ is stationary.
Similarly, we can find $C \subseteq M$ finite, 
such that $p_M$ and $q_M$ are based on $C$
and $a \fork_{CE} b$.
Since $C \subseteq M$ finite and $M \in \mathcal{K}$, 
there exists $a^*,b^*,E^* \subseteq M$,
such that $\tp(abE,C)=\tp(a^*b^*E^*,C)$, and so 
$a^* \fork_{CE^*} b^*$.
Since $\tp(ab, N)$ is based on $E$, then $\tp(ab,CE)$ is 
stationary based on $E$, 
so
$\tp(a^*b^*,CE^*)$ is stationary based on $E^*$.
Therefore, we can choose $a'b' \models \tp(a^*b^*,CE^*)_M$, and by choice 
of $C$, necessarily $a' \models p_M$ and $b' \models q_M$.

Hence, by assumption on $p_M, q_M$, we have $a' \nonfork_M b'$,
so also $a' \nonfork_{CE^*} b'$.
But this implies $a^* \nonfork_{CE^*} b^*$, by choice of $a'b'$, 
a contradiction.
\epf

\begin{lemma} \label{technical} Let $N$ be $(D,\mu)$-homogeneous.
If $a \nonfork_N b$ and $\tp(a, Nb)$ is $D^s_{\mu}$-isolated, then $a \in N$.
\end{lemma}

\bpf Since $p=\tp(a, Nb)$ is $D^s_{\mu}$-isolated, 
there is $C \subseteq N$, $|C| < \mu$ such that $\tp(a, Cb)$ isolates $p$.
Since $\tp(b, N)$ is stationary, we may assume that $\tp(b,N)$ 
does not split over $C$.
Since, by Theorem \ref{forkingproperties} \ref{symm} also $b \nonfork_N a$, so
we may assume that $\tp(b, Na)$ does not split over $C$. 

Now, since $N$ is $(D,\mu)$-homogeneous, there is $a' \in N$, such
that $\tp(a,C)=\tp(a', C)$. 
But since $\tp(b, Na)$ does not split over $C$, 
then $\tp(ab,C)=\tp(a'b,C)$.
Hence $\tp(a, N)=\tp(a',N)$, so that $a \in N$.
\epf

We recall a definition from \cite{Sh:1}.

\begin{definition} A $D$-model $M$ is 
\emph{maximally $(D,\mu)$-homogeneous} if
$M$ is $(D, \mu)$-homogeneous, but not $(D, \mu^+)$-homogeneous.
\end{definition} 

\begin{theorem} \label{somemodels}
Suppose $D$ is not unidimensional.
Then there is a
maximally $(D,\mu)$-homogeneous model $M$ of cardinality $\lambda$,
for every $\lambda \geq \mu \geq |T|+|D|$.
\end{theorem}
\bpf Suppose $D$ is totally transcendental and not unidimensional.
Then there exists $M,  N$ in  $\mathcal{K}$ and a minimal
type
$q(x, \bar{a})$ over $M$ with the property that
\begin{equation} \tag{*}
q(M, \bar{a}) = q(N, \bar{a})
\quad
\text{ and }
\quad
M \subseteq N,
\quad 
M \not = N.
\end{equation}

Using the Downward L\"{o}wenheim Skolem Theorem and prime models, 
we may assume that $|q(M,\bar{a})| \leq |T| + |D|$.
Let $\lambda \geq \mu \geq |T|+|D|$ be given.
We first show that we can find $M$, $N \in \mathcal{K}$
satisfying (*) such that in addition
$\|M\|=|q(M,\bar{a})| =\mu$.

Since $M \not =N \in \mathcal{K}$, there is $b \in N-M$, so
$p=\tp(b,M) \in S_D(M)$ is big and stationary.
This implies that $a' \nonfork_M b'$ for any $a' \models q_M$ 
and $b' \models p$
(by an automorphism sending $b'$ to $b$, 
it is enough to see $a' \nonfork_M b$, but
this is obvious, otherwise $\tp(a',Mb)$ is not big, 
thus cannot be big for $N$ by Lemma \ref{bigone},
hence it has to be realized in $N-M$, 
which implies that $a' \in N-M$, contradicting
$q(M,\bar{a})=q(N,\bar{a})$).

Construct $\langle M_i \mid i \leq \mu \rangle$ increasing and 
$I = \{ a_i \mid i < \mu \}$, $a_i \not \in M_i$
realizing $q_{M_i}$,
such that:
\begin{enumerate}
\item 
$M_{i+1} \in \mathcal{K}$ is $D^s_{\aleph_0}$-primary over $M_i \cup a_i$;
\item 
$M_0 = M$;
\item 
$M_i = \Union_{j<i} M_j$ when $j$ is a limit ordinal;
\item 
If $b'$ realizes $p_{M_i}$, and $N^*$ is 
$D^s_{\aleph_0}$-primary over $M_i \cup b'$,
then $q(M_i,\bar{a})=q(N^*,\bar{a})$.
\end{enumerate}

This is enough: Consider $N$ $D^s_{\aleph_0}$-primary over $M_\mu \cup b'$,
where $b' \models p_{M_\mu}$. 
Then $b' \in N-M_\mu$ and yet 
$q(M_\mu,\bar{a})=q(N,\bar{a})$, so (*) holds. 
Furthermore, $\|M_\mu\|=|q(M_\mu,\bar{a})|=\mu$.

This is possible: 
\begin{itemize}
\item
For $i=0$, this follows from the definition of $q$ 
(send $b'$ to $b$ by an automorphism, 
fixing $M$, to obtain a realization of $q_M$ in $N-M$).

\item
If $i$ is a limit ordinal, and $b' \models p_{M_i}$, then
this implies that $b' \models p_{M_j}$, for any $j < i$. 
Also, if $N^*$ is prime over $M_i \cup b'$, and $c \in N^* - M_i$
realizes $q(x, \bar{a})$, then $\tp(c, M_i b')$ is 
$D^s_{\aleph_0}$-isolated over some $\bar{m}b$, 
and $\bar{m}b \in M_j$ for some $j < i$, hence
$c \in M_j$ by induction hypothesis, a contradiction.

\item
For $i=j+1$. Let $b' \models p_{M_j}$ and $N^*$ be prime over $M_j \cup b'$.
Suppose that $c \in N^*-M_j$ realizes $q(x, \bar{a})$.
Then, since $c \not \in M_j$,  we must have $\tp(c, M_j)$ is big, so 
$c \models q_{M_j}$.
Hence, by Lemma \ref{tech} we have $c \nonfork_{M_j} b'$.
But $\tp(c, M_j b')$ is $D^s_{\aleph_0}$-isolated, 
so by Lemma \ref{technical}, we must have $c \in M_j$, a contradiction.
Hence $q(M_i)=q(N^*)$ and we are done.
\end{itemize}

Let $M^* = M_\mu$, and fix $b \models p_{M^*}$.
We now show that we can find a $(D,\mu)$-homogeneous
model $N \in \mathcal{K}$ of cardinality $\lambda$ such that
$M^*$ and $N$ satisfy (*).
This implies the conclusion of the theorem:  $N$ is
$(D,\mu)$-homogeneous of cardinality $\lambda$;
$N$ is \emph{not} $(D,\mu^+)$-homogeneous, since 
$N$ omits $q_{M^*} \in S_D(M^*)$, and $\| M^* \| = \mu$.

We construct $\langle N_i \mid i \leq \lambda \rangle$ increasing, and
$b_i \not \in N_i$ realizing $p_{N_i}$  such that:
\begin{enumerate}
\item 
$b_0=b$ and $N_0$ is $D^s_{\mu}$-primary over $M^* \cup b$;
\item 
$N_{i+1}$ is $D^s_{\mu}$-primary over $N_i \cup b_i$;
\item 
$N_i= \Union_{j<i} N_i$, when $i$ is a limit ordinal;
\item 
$\|N_i\|\leq \lambda$;
\item 
$N_i$ is $(D, \mu)$-homogeneous;
\item 
$q(N_i,\bar{a})=q(M^*,\bar{a})$.
\end{enumerate}

This is clearly enough: $N_\lambda$ is as required.

This is possible: We construct $N_i$ by induction on $i \leq \lambda$.

\begin{itemize}
\item
For $i=0$, let $N^*  \subseteq N_0$ be 
$D^s_{\aleph_0}$-primary over $M^* \cup b$. 
We have $q(N^*,\bar{a})=q(M^*,\bar{a})$ by construction of $M^*$, so
it is enough to show that $q(N^*,\bar{a})=q(N_0,\bar{a})$.
Suppose not and let $c \in N_0 - N^*$ realize $q(x, \bar{a})$.
Then, $c$ realizes $q_{N^*}$ since $\tp(c,N^*)$ is big, 
and further there is $A \subseteq M^*$,
$|A|<\mu$ such that $\tp(c, Ab)$ isolates $\tp(c,M^*b)$. 
By Lemma \ref{averagemorley} since $I$ is based on $q$, we
have $\av(I,N^*)=q_{N^*}$,
where $I=\{a_i \mid i<\mu \} \subseteq M^*$ defined above.
But since both $\tp(c, Ab)$ and $\tp(c,M^*)$ are big, we must have 
$\tp(c,Ab)=\av(I, Ab)$ and $\tp(c,M^*)=\av(I, M^*)$.
Hence $\av(I,Ab) \vdash \av(I,M^*)$.
Now, by Lemma \ref{kappa}, 
we can find $I' \subseteq I$, $|I'|<\mu$ such that  $I-I'$ is indiscernible
over $Ab$. 
Since $|I|=\mu$, then $I-I' \not = \empty$ and
 all elements of $I-I'$ realize $\av(I,Ab)$,
hence also $\av(I,M^*)=q_{M^*}$.
But this is impossible since $I \subseteq M^*$.
Therefore $q(N_0,\bar{a})=q(N^*,\bar{a})=q(M^*, \bar{a})$.
\item
For $i$ a limit ordinal, the only condition to check is that 
$N_i$ is $(D, \mu)$-homogeneous,
but this follows from Theorem \ref{transunion}.
\item
For $i=j+1$, by induction hypothesis, 
we have $q(N_j,\bar{a})=q(M^*, \bar{a})$,
so it is enough to show that $q(N_{j+1}, \bar{a})=q(N_j,\bar{a})$.
Suppose $c \in N_{j+1}$ realizes $q$.
Since $N_{j+1}$ is $D^s_{\mu}$-primary over $N_j\cup b_j$, 
we have $\tp(c, N_j\cup b_j)$ is $D^s_{\mu}$-isolated.
But $c \nonfork_{N_j} b_j$, by Lemma \ref{tech}.
Therefore, by Lemma \ref{technical}, we have that $c \in N_j$.
This shows that $q(N_{j+1},\bar{a})=q(M^*,\bar{a})$.
\end{itemize}

This completes the proof.
\epf

\begin{corollary} \label{catuni}
Let $D$ be totally transcendental.
If $\mathcal{K}$ is categorical in some
$\lambda > |T|+|D|$ then $D$ is unidimensional.
\end{corollary}
\bpf Otherwise, there is a $D$-homogeneous model of cardinality $\lambda$
and a maximally $(D, |T|+|D|)$-homogeneous model of cardinality $\lambda$.
Hence $\mathcal{K}$ is not categorical in $\lambda$, since these models
cannot be isomorphic.
\epf

We now obtain strong structural results when $D$ is unidimensional.

\begin{theorem} \label{minimal} 
Let $D$ be  unidimensional.
Then every $M \in \mathcal{K}$ is prime and minimal over $q(M,\bar{a})$,
for any minimal type $q(x, \bar{a})$ over $M$.
\end{theorem}
\bpf Let $M \in \mathcal{K}$ be given.
Since $D$ is totally transcendental, 
there exists a minimal type $q(x, \bar{a})$ over $M$.
Consider $A= q(M, \bar{a})$.
To check minimality, suppose there was $N \in \mathcal{K}$, such that
$A \subseteq N \subseteq M$. 
Since $q(N,\bar{a})=A=q(M,\bar{a})$,
we must have $N=M$, by unidimensionality of $D$. 
We now show that $M$ is prime over $A$. 
Since $D$ is totally transcendental, there is $M^* \in \mathcal{K}$ 
prime over $A$.
Hence, we may assume that $A \subseteq M^* \subseteq M$.
Now the minimality of $M$ implies that $M=M^*$, so $M$ is prime over $A$.
Clearly, any other minimal type would have the same property.
\epf

We next establish two lemmas, which are key results to carry out the
geometric argument for the categoricity theorem.

\begin{lemma} \label{saturation} Let 
 $M \in \mathcal{K}$ and suppose
that  $q(x, \bar{a})$ is minimal over $M$.
If  $W=q(M,\bar{a})$ has dimension $\lambda$ infinite,
then $W$ realizes every extension $p \in S_D(A)$ of 
type $q$, provided  $A$ is a subset of $W$ of cardinality
less than the dimension $\lambda$.
\end{lemma}

\bpf Let $p \in S_D(A)$ be given extending $q$. 
Let $c \in \mathfrak{C}$ realize $p$. 
If $p$ is not big for $M$, then $p$ is not realized outside $M$
so $c \in M$. Hence $c \in W$ since $p$ extends $q$.
If however $p$ is big for $M$, then $p$ is big and then  
by Lemma \ref{baseaverage} and Theorem \ref{minimalregular} we have
that 
$p=\av(I,A)$, where $I$ is any basis of $W$ of cardinality $\lambda$. 
But $|I| =\lambda \geq |A|^+ +\aleph_0$, so by Lemma \ref{kappa}
and definition of averages, $\av(I,A)$ is realized 
by some element of $I \subseteq W$.
Hence $p$ is realized in $W$.
\epf

\begin{lemma} \label{dimension} Let $D$ be unidimensional and 
let $M$ be in $\mathcal{K}$ of cardinality $\lambda > |T|+|D|$.
Suppose $q(x,\bar{a})$ is minimal over $M$.
Then $q(M, \bar{a})$ has dimension $\lambda$.
\end{lemma}  
\bpf Let $M\in \mathcal{K}$ be given and $q(x,\bar{a})$ be minimal.
Construct $\langle M_\alpha \mid \alpha < \lambda \rangle$ 
strictly increasing and continuous
such that $\bar{a} \in M_0$, $M_\alpha \subseteq M$ and 
$\|M_\alpha \| = |\alpha| + |T|+|D|$.

This is possible by Theorem \ref{prime}: 
For $\alpha=0$, just choose $M_0 \subseteq M$ prime over $\bar{a}$.
For $\alpha$ a limit ordinal, let $M_\alpha = \Union_{\beta<\alpha} M_\beta$.
At successor stage,
since $\|M_\alpha\| \leq |\alpha| + |T| +|D| < \lambda$, 
there exists $a_\alpha \in M -M_\alpha$,
so we can choose $M_{\alpha+1} \subseteq M$ prime over 
$M_{\alpha} \cup a_\alpha$.

This is enough: Since $D$ is unidimensional, 
we can find $c_{\alpha} \in M_{\alpha+1}-M_{\alpha}$
realizing $q$. 
By definition, 
$\tp(c_\alpha, \Union \{ c_\beta \mid \beta < \alpha \})$ is big, since
$c_\alpha \not \in M_\alpha$. 
Hence $c_\alpha \not \in cl( \Union \{ c_\beta \mid \beta < \alpha \})$.
Therefore $\{ c_\alpha \mid \alpha < \lambda\}$ is independent 
and so $q(M, \bar{a})$
has dimension at least $\lambda$. 
Hence since $\|M\|=\lambda$, then
$q(M,\bar{a})$ has dimension $\lambda$.
\epf

\begin{theorem} \label{unicat} 
Let $D$ be unidimensional.
Then $\mathcal{K}$ is categorical in every $\lambda  > |T|+|D|$.
\end{theorem} 
\bpf Let $M_l \in \mathcal{K}$ for $l=1,2$ be of cardinality
 $\lambda > |T|+|D|$.
Since $D$ is totally transcendental, we can choose, $q(x, \bar{a}_1)$ minimal, 
with $\bar{a}_1 \in M_1$. 
Now, since $M_2$ is $(D, \aleph_0)$-homogeneous, 
we can find $\bar{a}_2 \in M_2$ 
such that $\tp(\bar{a}_1, \empty)=\tp(\bar{a}_2, \empty)$.
Then $q(x,\bar{a}_1)$ is minimal also.
Let $W_l=q(M_l, \bar{a}_l)$ for $l=1,2$.
Since $D$ is unidimensional, by Lemma \ref{dimension}, 
we have $\dim(W_l)= \lambda > |T|+|D|$.
Hence, by Lemma \ref{saturation}
 every type extending $q(x, \bar{a}_l)$ over a subset of 
$W_l$ of cardinality less than $\lambda$
is realized in $W_l$, for $l=1,2$. 
This allows us to construct by induction an elementary mapping 
$g$ from $W_1$ onto
$W_2$ extending $\langle \bar{a}_1,\bar{a}_2\rangle$.
By Theorem \ref{minimal}, $M_l$ is prime and minimal over $W_l$, for $l=1,2$.
Hence, in particular $M_1$ is prime over $W_1$, so there is 
$f:M_1 \rightarrow M_2$ elementary extending $g$. 
But now $\rang( f) $ is a $(D, \aleph_0)$-homogeneous model
containing $W_2$, so by minimality of $M_2$ over $W_2$ we have
$\rang( f)= M_2$.
Hence $f$ is also onto, and so $M_1$ and $M_2$ are isomorphic.
\epf

We can now summarize our results.

\begin{corollary} \label{corollary1} 
Let $D$ be totally transcendental. The following conditions are equivalent:
\begin{enumerate}
\item 
$\mathcal{K}$ is categorical in every $\lambda > |T|+|D|$;
\item 
$\mathcal{K}$ is categorical in some $\lambda > |T|+|D|$;
\item 
$D$ is unidimensional;
\item 
Every $M \in \mathcal{K}$ is prime and minimal over $q(M,\bar{a})$,
where $q(x,\bar{a})$ is any minimal type over $M$;
\item 
Every model $M \in \mathcal{K}$ of cardinality $\lambda > |T|+|D|$
is $D$-homogeneous.
\end{enumerate}
\end{corollary}
\bpf  
\begin{nolabels}
\item
(1) implies (2) is trivial.
\item
(2) implies (3) is Theorem \ref{catuni}.
\item
(3) implies (1) is Theorem \ref{unicat}.
\item
(3) implies (4) is Theorem \ref{minimal}.
\item
(4) implies (3) is clear since prime models exist by Theorem \ref{prime}.
\item
(5) implies (1) is by back and forth construction, 
similarly to the corresponding proof with saturated models.
\item
(1) implies (5) since for each $\lambda > |D|+|T|$ there exist a 
$(D,\lambda)$-homogeneous model of cardinality $\lambda$ 
(e.g. by Theorem \ref{prime}).
\end{nolabels}
\epf

\begin{corollary}
Let $D$ be totally transcendental.
If $\mathcal{K}$ is not categorical in some $\lambda_1 > |T|+|D|$, then
\begin{enumerate}
\item
If $T$ is countable, then there are at least $|\alpha|$ models
of cardinality $\aleph_\alpha$ in $\mathcal{K}$;
\item
For every $\lambda \geq \mu \geq |T|+|D|$
there is a maximally $(D,\mu)$-homogeneous of cardinality $\lambda$.
\end{enumerate}
\end{corollary}
\bpf (1) follows from (2).
For (2), notice that $D$ is not unidimensional by 
above Corollary, so 
the result follows from Theorem \ref{somemodels}.
\epf

\begin{corollary} \label{corollary3} 
Let $D$ be totally transcendental.
Suppose there is a maximally $(D,\mu)$-homogeneous model 
of cardinality $\lambda > |T|+|D|$
for some $\lambda > \mu \geq \aleph_0$.
Then for every $\lambda \geq \mu \geq |T|+|D|$
there is a maximally $(D,\mu)$-homogeneous of cardinality $\lambda$.
\end{corollary}

\bpf Notice that $M \in \mathcal{K}$, and so $\mathcal{K}$ 
is not categorical in $\lambda$.
Hence, by the previous corollary, $D$ is not unidimensional, 
so the result follows
from Theorem \ref{somemodels}.
\epf

As a last Corollary, we obtain a generalization
of Keisler's Theorem (notice that $\mathcal{K}$
is the class of atomic models in this case).
 We \emph{do not} 
assume that $D$ is totally transcendental.

\begin{corollary} \label{corollary2} Let $|T| < 2^{\aleph_0}$,
and suppose $D$ is the set of isolated types of $T$.
The following conditions are equivalent.
\begin{enumerate}  
\item 
$\mathcal{K}$ is categorical in every $\lambda > |T|$;
\item 
$\mathcal{K}$ is categorical in some $\lambda > |T|$;
\item 
$D$ is totally transcendental and unidimensional;
\item
$D$ is totally transcendental and every model of $\mathcal{K}$
is prime and minimal over $q(M,\bar{a})$,
where $q(x,\bar{a})$ is any minimal type over $M$;
\item Every model $M \in \mathcal{K}$ of cardinality $\lambda > |T|+|D|$
is $D$-homogeneous.
\end{enumerate} 
\end{corollary}
\bpf (5) implies (1) and (2) by back and forth construction.
The rest of the proof follows from
\ref{corollary1}, since
conditions (1), (2), (3) and (4) imply
that $D$ is totally transcendental.
More precisely (1) and (2) imply that $D$ is 
stable in $|T| < 2^{\aleph_0}$ and
hence totally transcendental: this
is a standard fact using Ehrenfeucht-Mostowski models.
For (3) and (4) it is a hypothesis.
\epf

\end{document}